\documentclass[leqno]{article}
\usepackage{amssymb,amsfonts,amsmath,theorem,longtable,url}
\usepackage[a4paper,left=3cm,right=3cm,top=3cm,bottom=3cm]{geometry}

\newcommand{\FF}{\mathbb F}

\newcommand{\q}{\textrm{q}}

\newcommand{\cpave}{\overline{c_p}}

\newtheorem{lemma}{Lemma}
\newtheorem{corollary}{Corollary}
\newtheorem{theorem}{Theorem}

\newtheorem{conjecture}{Conjecture}
\newtheorem{principle}{Principle}
{\theorembodyfont{\rmfamily}}

\begin{document}

\title{The probability that the number of points on the Jacobian of a genus $2$ curve
is prime}
\author{Wouter Castryck, Amanda Folsom, Hendrik Hubrechts,\\ and
Andrew V.\ Sutherland}
\date{}

\maketitle

\begin{abstract}
In 2000,
Galbraith and McKee heuristically derived
a formula that estimates the probability
that a randomly chosen elliptic curve over a fixed finite prime field has
a prime number of rational points.
We show how their heuristics can be generalized to
Jacobians of curves of higher genus. We then elaborate this in genus $g=2$ and
study various related issues, such as the probability
of cyclicity and the probability of primality of the number of points on the curve itself. Finally, we discuss the asymptotic behavior for $g \rightarrow \infty$. \\

\noindent MSC 2010: 11N05, 11G10, 11G20
\end{abstract}

\section{Introduction and overview}

\subsection{The Galbraith-McKee conjecture: elliptic curves} \label{intro1}

In \cite{GalbraithMcKee}, Galbraith and McKee studied the
probability that a randomly chosen elliptic curve over a finite prime field has
a prime number of rational points. They conjectured the following. For a prime number $p > 3$, let $P_1(p)$ be
the probability that a uniformly randomly chosen integer in the Hasse interval $[p + 1 - 2\sqrt{p}, p + 1 + 2\sqrt{p} ]$
is prime. Let $P_2(p)$ be the probability that the elliptic curve defined by $y^2 = x^3 + Ax + B$, for a
uniformly randomly chosen pair $(A,B)$ in the set
\[ \mathcal{H}_{AB} = \left\{ \left. \, (A,B) \in \mathbb{F}_p^2 \, \right| \, 4A^3 + 27B^2 \neq 0 \, \right\}, \]
has a prime number of rational points (including the point at infinity).
\begin{conjecture}[Galbraith-McKee \textbf{\cite[Conjecture~A]{GalbraithMcKee}}] \label{galbraithmckee} Define
\[ c_p = \frac{2}{3} \cdot \prod_{\ell > 2} \left( 1 - \frac{1}{(\ell - 1)^2} \right) \cdot \prod_{\ell \mid p-1, \ell > 2} \left( 1 + \frac{1}{(\ell + 1)(\ell -2)} \right), \]
where the products are over all primes $\ell$ satisfying the stated conditions. Then
\[ \lim_{p \rightarrow \infty} \left( P_2(p)/P_1(p) - c_p \right) = 0. \]
\end{conjecture}
The constant $c_p$ lies
between 0.44010 and 0.61514. In general, the conjecture predicts that
elliptic curves are about half as likely to have prime orders as one might expect.

The study of the probability of primality is partly motivated by
elliptic curve cryptography. For an elliptic curve over a finite
field to be suitable as the underlying group for
Diffie-Hellman key exchange, its number of rational points is
preferably prime (although small cofactors are often tolerated). In practice, a `good' elliptic
curve is often found by repeatedly counting the number of rational
points on randomly chosen elliptic curves, e.g.\ using the SEA
algorithm \cite{SEA}, until a prime number is hit. The above conjecture
predicts that this process works slightly worse than one would
naively assume.

Galbraith and McKee provided both experimental support and heuristic evidence in favor
of Conjecture~\ref{galbraithmckee}. Their main argument uses the Hurwitz-Kronecker class number formula,
which counts bivariate quadratic forms up to equivalence.
A second argument estimates the probability of primality
by naively multiplying the expected probabilities of being coprime to $2, 3, 5, 7, 11, \dots$
For elliptic curve orders, these expected probabilities were devised by Lenstra \cite[Proposition~1.14]{lenstra}.
When taking the quotient of the resulting estimates for $P_2(p)$ and $P_1(p)$, one exactly finds $c_p$.
A reasoning of this kind had already been made by Koblitz \cite[p.\ 160]{Koblitz}
in the dual setting where one fixes an elliptic curve over $\mathbb{Q}$ and reduces it modulo varying primes --- a similar discussion
on the case where one fixes a CM-curve of genus $2$ over $\mathbb{Q}$ can be read in Weng's thesis \cite[Section~5.2]{WengThesis}.
Galbraith and McKee called their second heuristics `not very honest', however, due
to subtleties reflected in Mertens' theorem. We will discuss these
in Section~\ref{framework}.

\subsection{Genus 2 curves} \label{intro2}

Nonetheless -- and this may be thought of as an underlying meta-conjecture --
these second heuristics work very well in practice, as is confirmed
experimentally in Section~\ref{experimental}. Moreover,
they seem more flexible towards generalizing
Conjecture~\ref{galbraithmckee} to Jacobians of curves of higher genus,
which have also been proposed for use in cryptography.
The required analogues of Lenstra's theorem are provided
by a recursive formula due to Achter and Holden \cite[Lemma~3.2]{AchterHolden},
which we turn into a closed expression in Section~\ref{countingmatrices}.

In this article, we elaborate this for curves of genus $2$, which is the most relevant case for cryptography.
We derive the following conjecture.
For a prime number $p > 2$, let $P_1(p)$ be the probability that a
uniformly randomly chosen
integer in the Hasse-Weil interval
\[ \left[(\sqrt{p} - 1)^4, (\sqrt{p} + 1)^4 \right] \]
is prime. Let $P_2(p)$ be the probability that the Jacobian of
the genus $2$ curve defined by $y^2 = f(x)$, for a randomly chosen polynomial $f(x)$ in the set
\[ \mathcal{H}_6 = \left\{ \left. \, f(x) \in \mathbb{F}_p[x] \, \right| \, f(x) \text{ square-free of degree $6$} \, \right\}, \]
has a prime number of rational points.
\begin{conjecture}[see Section~\ref{lenstrageneralization}] \label{ourconjecture} Define
\[ c_p = \frac{38}{45} \ \cdot \ \prod_{\ell > 2} \left( 1 - \frac{\ell^2 - \ell - 1}{(\ell^2-1)(\ell - 1)^2} \right) \ \cdot \ \prod_{\ell \mid p-1, \ell > 2}
\left( 1 + \frac{\ell^4 - \ell^3 - \ell  - 2}{(\ell^3 - 2\ell^2 - \ell +3)(\ell^2 + 1)(\ell + 1)} \right) , \]
where the products are over all primes $\ell$ satisfying the stated conditions. Then
\[ \lim_{p \rightarrow \infty} \left( P_2(p)/P_1(p) - c_p \right) = 0. \]
\end{conjecture}
We implicitly assume that $P_1(p) \neq 0$ for all $p$,
which is an open problem in its own (see \cite[Section~2.2]{Broker} for a related discussion). The constant $c_p$ lies between 0.63987 and
0.79890.
Summarizing,
in genus 2 prime order Jacobians are also slightly disfavored, but to
a lesser extent than in genus 1.

\subsection{Averaging over $p$} \label{averagingoverp}

By averaging $c_p$ over all primes $p$, it becomes meaningful to measure
the prime-disfavoring behavior by a single constant. For elliptic curves, this gives:

\begin{lemma}[see Section~\ref{framework}] \label{galbraithmckeeaverage}
For each prime $p > 3$, let $c_p$ be as in Conjecture~\ref{galbraithmckee}.
Then
\[\cpave = \lim_{n \rightarrow \infty} \frac{1}{\pi(n)} \sum_{3 < p \leq n} c_p = \prod_{\ell} \left( 1 - \frac{\ell^2 - \ell - 1}{(\ell^2 - 1)(\ell - 1)^2 } \right) \approx 0.50517.\]
Here, $\pi$ is the prime-counting function, and the product is over all primes $\ell$.
\end{lemma}
This confirms a constant
obtained by Koblitz \cite[p.~160]{Koblitz} and subsequently verified by Balog, Cojocaru and David \cite[Theorem~1]{BCD}.
In genus $2$, the average reads:

\begin{lemma}[see Section~\ref{lenstrageneralization}] \label{ourconjectureaverage}
For each prime $p > 2$, let $c_p$ be as in Conjecture~\ref{ourconjecture}. Then
\[\cpave = \lim_{n \rightarrow \infty} \frac{1}{\pi(n)} \sum_{2 < p \leq n} c_p = \prod_{\ell} \left( 1 - \frac{\ell^6 - 2\ell^5 + 3\ell + 1}{(\ell^2 - 1)^2 (\ell^2 + 1)(\ell - 1)^2} \right) \approx 0.69464,\]
where again the product is over all primes $\ell$.
\end{lemma}

\subsection{Imposing a rational Weierstrass point}

Instead of using $\mathcal{H}_6$, we can
choose $f(x)$ uniformly at random from the set
\[ \mathcal{H}_5^\text{m} = \left\{ \left. \, f(x) \in \mathbb{F}_p[x] \, \right| \, f(x) \text{ monic and square-free of degree $5$} \, \right\}. \]
This situation matches better with common cryptographic practice.
However, it alters the notion of taking a random genus 2 curve,
since here one imposes the existence of a rational Weierstrass point.
As before, for each prime $p>2$, let $P_1(p)$ be the probability that
a uniformly randomly chosen integer in the Hasse-Weil interval
$[(\sqrt{p} - 1)^4, (\sqrt{p} + 1)^4]$ is prime,
but now let $P_2(p)$ be the probability that a random genus $2$ curve, in the above sense,
has a Jacobian with a prime number of rational points.

\begin{conjecture}[see Section~\ref{rationalweierstrass}] \label{imposingweierstrass}
Let $c_p$ be as in Conjecture~\ref{ourconjecture}. Then
\[ \lim_{p \rightarrow \infty} \left( P_2(p)/P_1(p) - \frac{9}{19}c_p \right) = 0. \]
\end{conjecture}
The constant $\frac{9}{19}c_p$ lies between $0.30309$ and $0.37843$, so prime orders
become dramatically less probable. This is entirely due to the fact that
the probability of having rational $2$-torsion increases from $\frac{26}{45}$ to $\frac{4}{5}$.
In Section~\ref{rationalweierstrass}, we will illustrate why for odd $\ell$, the expected probability
of having rational $\ell$-torsion is most likely unaffected.

Averaging $\frac{9}{19}c_p$ over all primes $p$ as
in Section~\ref{averagingoverp} gives approximately $0.32904$ (i.e.\ $\frac{9}{19}$ times the constant of Lemma~\ref{galbraithmckeeaverage}).

\subsection{The number of points on the curve itself} \label{curveitselfintro}

We can also estimate the probability that
the number of rational points on the curve itself, rather
than its Jacobian, is prime.
For each prime $p>2$ and with $f(x)$ chosen uniformly at random
from $\mathcal{H}_6$,
let $P_2(p)$ be the probability that the nonsingular complete model of $y^2 = f(x)$ has
a prime number of rational points. Let $P_1(p)$ be the probability
that an integer, chosen uniformly at random from the Hasse-Weil interval
\[ [p+1 - 4\sqrt{p}, p+1 + 4\sqrt{p}] \]
is prime. For $\ell \neq p$ prime, define
\begin{eqnarray*}
a_{\ell,p} &:=& \#\{(x,y)\in\mathbb{F}_\ell^\times\times(\mathbb{F}_\ell^\times\backslash \{-p\}) \ |\ (x+y/x)(1+p/y)=p+1\},\\
\beta_{\ell,p} &:=& (\ell-1)(\ell^5 - \ell^3 + 2) - a_{\ell,p} -\begin{cases}(\ell^3-1) & \textnormal{if $p\equiv -1\bmod \ell$,}\\ 0  & \textnormal{otherwise.}\end{cases}
\end{eqnarray*}

\begin{conjecture}[see Section~\ref{curveitself}] \label{curveitselfconj}
Define
\[ c_p = \frac{38}{45} \prod_{\substack{\ell>2 \\ \ell \neq p \\}} \frac{\ell \cdot \beta_{\ell,p} }{(\ell^4-1)(\ell^2-1)(\ell-1) } \]
where the product is over all primes $\ell>2$. Then
\[ \lim_{p \rightarrow \infty} \left( P_2(p)/P_1(p) - c_p \right) = 0. \]
\end{conjecture}
The constant $c_p$ lies between $0.79605$ and $0.86548$, with an estimated average (in the sense of
Section~\ref{averagingoverp}) of $\overline{c_p} \approx 0.83376$. When switching to $\mathcal{H}_5^\text{m}$ instead of $\mathcal{H}_6$,
the leading factor $\frac{38}{45}$ should be replaced by $\frac{16}{15}$. The resulting constant
$c_p$ lies between $1.00553$ and $1.09323$, with an estimated average of $\overline{c_p} \approx 1.05317$, so prime orders actually become slightly favored.

\subsection{The probability of cyclicity}

Using similar heuristics, one can estimate for each prime $p>2$ the probability $P(p,2)$ that the group of
rational points on the Jacobian
of the curve defined by $y^2 = f(x)$, with $f(x)$ chosen uniformly at random from
$\mathcal{H}_6$,
is cyclic. This is done by considering for each prime $\ell$ the corresponding probability
for the $\ell$-torsion subgroup, and then taking the product.

For elliptic curves, one recovers a formula that was proven
by Vl\u{a}du\c{t}. Let $P(p,1)$ be the probability that the group of rational points on
a randomly chosen elliptic curve over $\mathbb{F}_p$ (as in Section~\ref{intro1})
is cyclic. Then:
\begin{theorem}[Vl\u{a}du\c{t} \textbf{\cite[Theorem~6.1]{vladut}}] \label{theoremvladut}
For each prime $p$, define
\[ c_p = \prod_{\ell \mid p-1} \left( 1 - \frac{1}{\ell(\ell^2-1)} \right)\]
where the product is over all primes $\ell$ satisfying $\ell \mid p-1$. Then
\[ \lim_{p \rightarrow \infty} (P(p,1) - c_p ) = 0.\]
\end{theorem}
The constant $c_p$ is contained in $[0.78816, 0.83334]$, with an average (in the sense of Section~\ref{averagingoverp}) of $\overline{c_p} \approx 0.81375$. In genus $2$, the same
reasoning gives:
\begin{conjecture}[see Section~\ref{cyclicity}] \label{cyclicitygenus2}
For each prime $p$, define
\[ c_p = \frac{151}{180} \cdot \prod_{\ell > 2, \ell \nmid p-1} \left( 1 - \frac{1}{\ell(\ell^2-1)(\ell-1)} \right) \cdot \prod_{\ell > 2, \ell \mid p-1}
\frac{\ell^8 - \ell^6 - \ell^5 - \ell^4  + \ell^2 + \ell + 1}{\ell^2(\ell^4-1)(\ell^2-1)},\]
where the products are over all primes $\ell$ satisfying the stated conditions. Then
\[ \lim_{p \rightarrow \infty} (P(p,2) - c_p ) = 0.\]
\end{conjecture}
The constant $c_p$ is contained in the interval $[0.79356,0.81918]$, with
average value $\cpave \approx 0.80883$. If we replace
$\mathcal{H}_6$ by $\mathcal{H}_5^\text{m}$, then the leading factor should be replaced
by $\frac{37}{60}$, in which case the constant $c_p$ is contained between $0.58335$ and $0.60218$, with average value
$\cpave \approx 0.59457$.

\subsection{Extension fields}

Fix a prime number $p$. Consider the alternative setup of finite fields $\mathbb{F}_{p^k}$ of growing extension degree $k$
over $\mathbb{F}_p$.
For $g \in \{1,2\}$, let $P_1(k,g)$ be the probability that a uniformly randomly chosen integer in the Hasse interval
$[(\sqrt{p^k}-1)^{2g}, (\sqrt{p^k} + 1)^{2g} ]$ is prime. Let $P_2(k,g)$ be the probability that
the Jacobian of the (hyper)elliptic curve defined by $y^2 + h(x)y = f(x)$, where the pair
$(h,f)$ is chosen from
\[ \begin{array}{c} \mathcal{H}_{g+1,2g+2} = \{ (f,h) \in \mathbb{F}_{p^k}[x] \times \mathbb{F}_{p^k}[x] \, | \, \deg h \leq g+1, \deg f = 2g+2, \\
\qquad \qquad \qquad \qquad \qquad \qquad \qquad y^2 + h(x)y = f(x) \text{ has geometric genus $g$}\} \\ \end{array} \]
uniformly at random, has a prime number of $\mathbb{F}_{p^k}$-rational points.

Then:

\begin{conjecture}[see Section~\ref{extensionfields}] \label{ext1}
Let
\[ c_k = \mu_p \cdot \prod_{\ell > 2} \left( 1 - \frac{1}{(\ell - 1)^2} \right) \cdot \prod_{\ell \mid p^k-1, \ell > 2} \left( 1 + \frac{1}{(\ell + 1)(\ell -2)} \right), \]
where the products are over all primes $\ell$ satisfying the stated conditions, and $\mu_p = 0$ if
$p = 2$ versus $\mu_p = \frac{2}{3}$ if $p>2$. Then
\[ \lim_{k \rightarrow \infty} \left( P_2(k,1)/P_1(k,1) - c_k \right) = 0. \]
\end{conjecture}
If $p>2$, the formula for $c_k$ closely matches the formula from $c_p$ from
Conjecture~\ref{galbraithmckee}, with $p^k-1$ in place of $p-1$, and takes values between $0.44010$ and $0.61514$. For $p=2$ we have $c_k=0$.
In genus $2$, the estimate reads:
\begin{conjecture}[see Section~\ref{extensionfields}] \label{ext2} Let
\[ c_k = \mu_p \ \cdot \ \prod_{\ell > 2} \left( 1 - \frac{\ell^2 - \ell - 1}{(\ell^2-1)(\ell - 1)^2} \right) \ \cdot \ \prod_{\ell \mid p^k-1, \ell > 2}
\left( 1 + \frac{\ell^4 - \ell^3 - \ell  - 2}{(\ell^3 - 2\ell^2 - \ell +3)(\ell^2 + 1)(\ell + 1)} \right) , \]
where the products are over all primes $\ell$ satisfying the stated conditions, and
$\mu_p = \frac{2}{3}$ if $p=2$ versus $\mu_p = \frac{38}{45}$ if $p > 2$. Then
\[ \lim_{k \rightarrow \infty} \left( P_2(k,2)/P_1(k,2) - c_k \right) = 0. \]
\end{conjecture}
Again for $p> 2$, the formula for $c_k$ matches the formula for $c_p$ in Conjecture~\ref{ourconjecture}
and takes values between $0.63987$ and $0.79890$. If $p=2$ then $c_k$ lies between $0.50516$ and
$0.63071$.

It is possible to average the above over $k$, where the result will depend on
the multiplicative orders of $p$ modulo the various $\ell$. Also, one can adapt Conjectures~\ref{ext1}
and~\ref{ext2}, and in fact any of the conjectures stated above, to the mixed case of just considering finite fields
$\mathbb{F}_q$ of growing cardinality.

\subsection{Asymptotics for growing genus} \label{asymptotics}

Instead of elaborating similar, increasingly complicated formulas for higher genera $g$,
we end with an analysis of the asymptotic behavior for $g \rightarrow \infty$. This may be of interest to
people studying analogues of the Cohen-Lenstra heuristics \cite{CohenLenstra,lengler} in the case of function fields,
though we will not push this connection. Note that due to computational limitations, the conjectures
below are no longer supported by experimental evidence and rely purely on the conjectured validity of our heuristic derivation.

For every prime number $p > 2$ and every integer $g \geq 1$, let $P_1(p,g)$ be the probability that a
uniformly randomly chosen
integer in the Hasse-Weil interval
\[ \left[(\sqrt{p} - 1)^{2g}, (\sqrt{p} + 1)^{2g} \right] \]
is prime. Let $P_2(p,g)$ be the probability that the Jacobian of
the genus $g$ curve defined by $y^2 = f(x)$, for a randomly chosen polynomial $f(x)$ in the set
\[ \mathcal{H}_{2g+2} = \left\{ \left. \, f(x) \in \mathbb{F}_p[x] \, \right| \, f(x) \text{ square-free of degree $2g+2$} \, \right\}, \]
has a prime number of rational points.

Then we have:
\begin{theorem}[see Section~\ref{lenstrageneralization}] \label{limithyperelliptic}
$\lim_{p, g \rightarrow \infty} P_2(p,g) = 0$.
\end{theorem}
Theorem~\ref{limithyperelliptic} holds
because the probability of having
rational $2$-torsion tends to $1$ as $g \rightarrow \infty$. However, this is a hyperelliptic phenomenon.
The limiting behavior becomes more interesting if instead one
defines $P_2(p,g)$ as the probability that the Jacobian of a random genus $g$ curve over $\mathbb{F}_p$ (e.g.
chosen from the set
\[ \mathcal{M}_g = \left\{ \, \text{curves of genus $g$ over $\mathbb{F}_p$} \, \right\} / \cong_{\mathbb{F}_p} \]
uniformly at random --- note that $\mathcal{M}_g$ is typically not well-understood) has a prime number of rational points. In this case, we expect:

\begin{conjecture}[see Section~\ref{lenstrageneralization}] \label{limitconjecture}
Define
\[ c_p = \frac{1}{\prod_{j=2}^\infty \zeta(j)} \cdot \prod_{\ell \mid p-1} \prod_{j=1}^\infty \frac{\ell^{2j}}{\ell^{2j} - 1}\]
where $\zeta$ is Riemann's zeta function and the product is over all primes $\ell$ satisfying the stated condition. Then
\[ \lim_{p,g \rightarrow \infty} \left( P_2(p,g)/P_1(p,g) - c_p \right) = 0. \]
\end{conjecture}
Again, we implicitly assume that $P_1(p,g)$ is nowhere zero. The constant $c_p$ lies
in the interval
\[ \left[ \frac{\prod_{j=1}^\infty \frac{2^{2j}}{2^{2j} - 1}}{\prod_{j=2}^\infty \zeta(j)},
 \frac{1}{\prod_{j=1}^\infty \zeta(2j+1)}
 \right] \subset [0.63287, 0.79353].\]
In other words, the prime-disfavoring effect persists as the genus grows. It even becomes
slightly more manifest than in genus $2$. A more
detailed analysis shows that the effect alternatingly strengthens and weakens as the genus
becomes odd and even, respectively.
As in Section~\ref{averagingoverp}, one can average $c_p$ over all primes $p>2$, yielding
a constant $\overline{c_p} \approx 0.68857$.

%Averaging $c_p$ over all primes $p > 2$ gives:
%\begin{conjecture}[see Section~\ref{lenstrageneralization}] \label{limitaverage}
%For each prime $p>2$, let $c_p$ be as in Conjecture~\ref{limitconjecture}. Then
%\[ \lim_{n \rightarrow \infty} \frac{1}{\pi(n)} \sum_{2 < p \leq n} c_p = \frac{1}{\prod_{j=2}^\infty \zeta(j)} \cdot
%\prod_\ell \frac{\ell - 2 + \prod_{j=1}^\infty \frac{\ell^{2j}}{\ell^{2j} - 1} }{\ell - 1} \approx 0.68857,  \]
%where the product is over all primes $\ell$.
%\end{conjecture}

Similarly, for every prime number $p>2$ and every integer $g \geq 1$, let $P(p,g)$
be the probability that the rational points of the Jacobian of the (hyper)elliptic curve $y^2=f(x)$,
with $f(x)$ picked from $\mathcal{H}_{2g+2}$ uniformly at random, constitute a cyclic group.

Then:
\begin{theorem}[see Section~\ref{cyclicity}] \label{limitcyclichyperelliptic}
 $\lim_{p,g \rightarrow \infty} P(p,g) = 0$.
\end{theorem}
Again, this is a hyperelliptic phenomenon due to $2$-torsion issues. If instead
we define $P(p,g)$ to be the probability that a curve chosen from $\mathcal{M}_g$ uniformly at random has
a cyclic Jacobian, then we expect:
\begin{conjecture}[see Section~\ref{cyclicity}] \label{limitcyclic}
Define
\[ c_p = \frac{1}{\prod_{j=2}^\infty \zeta(j)} \cdot \prod_{\ell \mid p-1} \prod_{j=1}^\infty \frac{\ell^{2j}}{\ell^{2j} - 1} \cdot \prod_{\ell \nmid p-1} \left( 1 + \frac{1}{\ell(\ell-1)} \right) \]
where $\zeta$ is Riemann's zeta function and the product is over all primes $\ell$ satisfying the stated conditions. Then
\[ \lim_{p,g \rightarrow \infty} \left( P(p,g) - c_p \right) = 0. \]
\end{conjecture}
Now the constant $c_p$ lies in the interval
\[ \left[  \frac{1}{\prod_{j=1}^\infty \zeta(2j+1)}, \frac{\prod_{j=1}^\infty \frac{2^{2j}}{2^{2j} - 1} \cdot \prod_{\ell >2} \left( 1 + \frac{1}{\ell(\ell -1)} \right)}{\prod_{j=2}^\infty \zeta(j)}
 \right] \subset [0.79352, 0.82004],\]
with an average (in the sense of Section~\ref{averagingoverp}) of $\overline{c_p} \approx 0.80924$.

%\subsection{Remarks} \label{remarks}

%The most important question that we do not address in this article, is
%how to adapt the foregoing to the case of primality \emph{up to a given cofactor $k$}.
%For this, one would need to generalize Lenstra's theorem to non-prime moduli. In the elliptic curve case, this was
%done by Howe \cite{Howe}, allowing Galbraith and McKee to adapt
%Conjecture~\ref{galbraithmckee} appropriately \cite[Conjecture~B]{GalbraithMcKee}. For
%higher genera, such a generalization is not available yet. It could in principle be achieved
%by counting matrices, using the random matrix model as described in Section~\ref{randommatrixgenus2}.
%But since the count does not seem straightforward, we leave this as a future research project.

\section{Common notions of randomness} \label{randomness}

By a randomly chosen (hyper)elliptic curve of genus $g \geq 1$ over a finite field $\mathbb{F}_q$ of odd characteristic, we
will usually mean the nonsingular complete model of a curve $y^2 = f(x)$, where $f$ is chosen from
\[ \mathcal{H}_{2g+2} = \left\{ \left. \, f(x) \in \mathbb{F}_q[x] \, \right| \, \text{$f(x)$ is square-free and $\deg f = 2g+2$} \, \right\} \]
uniformly at random.

Alternatively, one could take
the curve uniformly at random from
\[ \mathcal{M}_g^\text{hyp} = \{ \text{(hyper)elliptic genus $g$ curves over $\mathbb{F}_q$} \} / \cong_{\mathbb{F}_q}. \]
This randomness notion may be preferred from a theoretical point of view. It is fundamentally different from our first, in the sense that the map
\[ \mathcal{H}_{2g+2} \rightarrow \mathcal{M}_g^\text{hyp} : f \mapsto [y^2 = f(x)] \]
is not uniform. For small $q$ it does not even need to be surjective.
Therefore, the probability of having a certain geometric property may
change when moving from the one notion to the other. However,
as $q$ gets bigger and bigger, the change becomes negligible.
More precisely, for $q \rightarrow \infty$ ($g$ fixed), the proportion
of elements of $\mathcal{M}_g^\text{hyp}$ having $q(q^2-1)(q-1)/2$ pre-images in $\mathcal{H}_{2g+2}$ tends to $1$. This can
be elaborated following \cite[Section~1]{Nart}. Note that, despite the availability of a complete classification of (hyper)elliptic curves up to
$\mathbb{F}_q$-isomorphism \cite[Section~2]{Nart}, the set $\mathcal{M}_g^\text{hyp}$ is quite cumbersome
to work with.

Another setup, which is e.g.\ used in \cite[Theorem~3.1]{achter}, is to take $f$ uniformly at random from
\[ \mathcal{H}_{2g+2}^\text{m} = \left\{ \left. \, f(x) \in \mathbb{F}_q[x] \, \right| \, \text{$f(x)$ is monic, square-free and $\deg f = 2g+2$} \right\}, \]
instead of $\mathcal{H}_{2g+2}$. Again, this is different from either of the above notions. For small $q$, there may
exist curves having a model in $\mathcal{H}_{2g+2}$ that do not
have a model in $\mathcal{H}_{2g+2}^\text{m}$.
But again, as $q \rightarrow \infty$ ($g$ fixed), the difference dissolves. Indeed, consider the set
\[ \mathcal{S}_{2g+2} = \left\{ \left. \, (f, \alpha, \beta) \in \mathcal{H}_{2g+2} \times \mathbb{F}_q \times \mathbb{F}_q^\times \, \right| \,  f(\alpha) = \beta^2 \, \right\}.\]
Then we have a map
\[ \mathcal{S}_{2g+2} \rightarrow \mathcal{H}_{2g+2}^\text{m} : (f,\alpha,\beta) \mapsto \beta^{-2}x^{2g+2}f(1/x + \alpha),\]
which respects the isomorphism class of the corresponding curve, and which is onto and $q(q-1)$-to-$1$.
Therefore, taking $f$ uniformly at random from $\mathcal{H}_{2g+2}^\text{m}$ and using the $f$ of a uniformly
randomly chosen $(f,\alpha, \beta) \in S_{2g+2}$ give rise to equivalent randomness notions.
On the other hand, the map
\[ \mathcal{S}_{2g+2} \rightarrow \mathcal{H}_{2g+2} : (f,\alpha,\beta) \mapsto f \]
is asymptotically uniform, since
every $f \in \mathcal{H}_{2g+2}$ will have $q + O(\sqrt{q})$ pre-images by the
Hasse-Weil bound. This proves the claim.

In Section~\ref{extensionfields} we will allow $\text{char} \, \mathbb{F}_q = 2$ and use curves of the form $y^2 + h(x)y = f(x)$
with $(f,h)$ chosen from
\[ \begin{array}{c} \mathcal{H}_{g+1,2g+2} = \{ (f,h) \in \mathbb{F}_q[x] \times \mathbb{F}_q[x] \, | \, \deg h \leq g+1, \deg f = 2g+2, \\
\qquad \qquad \qquad \qquad \qquad \qquad y^2 + h(x)y = f(x) \text{ has geometric genus $g$}\} \\ \end{array} \]
uniformly at random. Again, it is easy to show that if $2 \nmid q$, the completing-the-square map
$\mathcal{H}_{g+1,2g+2} \rightarrow \mathcal{H}_{2g+2}$ is essentially uniform.

In this article, we will always consider statistical behavior for $q \rightarrow \infty$. In particular,
the validity of all statements below involving randomly chosen
curves in the sense of $\mathcal{H}_{2g+2}$ is preserved when switching to either of the above alternatives, and vice versa.
Some statements involve error terms, so in fact a more careful analysis is needed; we omit the details.

The picture does alter, however, when one takes
$f$ uniformly at random from
\[ \mathcal{H}_{2g+1} = \left\{ \, f \in \mathbb{F}_q[x] \, | \, \text{$f(x)$ is square-free and} \, \deg f = 2g+1 \, \right\}. \]
While this setting is often preferred in practice, this influences the story
as soon as $g \geq 2$, since it induces the existence of a rational Weierstrass point.
We will study this effect in detail for $g=2$ in Section~\ref{rationalweierstrass}.
On the other hand, writing
\[ \mathcal{H}_{2g+1}^\text{m} = \left\{ \, f \in \mathbb{F}_q[x] \, | \, \text{$f(x)$ is monic, square-free and} \, \deg f = 2g+1 \, \right\}, \]
the geometry-preserving map
\[ \mathcal{H}_{2g+1} \rightarrow \mathcal{H}_{2g+1}^\text{m}: f \mapsto \alpha^{2g} f(x/\alpha) \quad \text{(where $\alpha = \text{lc}(f)$)} \]
is onto and $(q-1)$-to-$1$. Therefore, $\mathcal{H}_{2g+1}$ and $\mathcal{H}_{2g+1}^\text{m}$ can be interchanged in
any probability statement below. If $g=1$ and moreover $3 \nmid q$, this also accounts
for
\[ \mathcal{H}_{AB} = \left\{ \left. \, (A,B) \in \mathbb{F}_q^2 \, \right| \, 4A^3 + 27B^2 \neq 0 \, \right\}, \]
since the completing-the-cube map $\mathcal{H}_{3} \rightarrow \mathcal{H}_{AB}$ is uniform.

Note that the sets $\mathcal{H}_{2g+2}, \mathcal{H}_{2g+2}^\text{m}, \mathcal{H}_{g+1,2g+2}, \mathcal{H}_{2g+1}, \mathcal{H}_{2g+1}^\text{m}, \mathcal{M}_g^\text{hyp}, \mathcal{H}_{AB}$ depend on $q$, while this is not included in the notation for sake of readability. However, it will always
be clear from the context which $q$ is used (it will typically be the prime number $p$ under consideration).

\section{Heuristic framework} \label{framework}

For prime numbers
$p >3$ and $\ell \neq p$, let $P(p, \ell)$
be the probability that the elliptic curve $E_{AB}$ defined by
$y^2 = x^3 + Ax + B$, for a randomly chosen pair $(A,B)$ in the
set $\mathcal{H}_{AB}$,
has $\ell$ dividing its number of rational points (including the point at infinity).
\begin{theorem}[Lenstra] \label{Lenstra}
There exist $C_1, C_2 \in \mathbb{R}_{>0}$, such that
\[ \left| P(p,\ell) - \frac{\ell}{\ell^2 - 1} \right| \leq C_1 \ell / \sqrt{p} \quad \text{if $\ell \mid p-1$} \quad \text{and} \]
\[ \left| P(p,\ell) - \frac{1}{\ell - 1} \right| \leq C_2 \ell / \sqrt{p} \quad \text{if $\ell \nmid p-1$} \]
for all pairs of distinct primes $p, \ell$ with $p > 3$.
\end{theorem}
\noindent \textsc{Proof.} See \cite[Proposition~1.14]{lenstra}, to which we refer for explicit estimates of the $C_i$. \hfill $\blacksquare$\\

We can now describe and discuss in more detail Galbraith and McKee's second heuristic argument supporting Conjecture~\ref{galbraithmckee}.
This is the type of reasoning behind all of our conjectures.
Let $\ell(p)$ be the largest prime for which $\ell(p) \leq \sqrt{p} + 1$. Let $n$ be an integer chosen uniformly at random
from the Hasse interval, and let $\eta$ be $\#E_{AB}(\mathbb{F}_p)$.
The aim is to estimate the ratio $P_2(p)/P_1(p)$, where $P_1(p)$ and $P_2(p)$
are as in Section~\ref{intro1}.
It can be rewritten as
\[ \frac{ P(2 \nmid \eta \text{ and } 3 \nmid \eta \text{ and } 5 \nmid \eta \text{ and } \dots \text{ and } \ell(p) \nmid \eta )}{P(2 \nmid n \text{ and } 3 \nmid n \text{ and } 5 \nmid n \text{ and } \dots \text{ and } \ell(p) \nmid n )}. \]
A first heuristic step is to approximate the above by
\[ \frac{P(2 \nmid \eta) P(3 \nmid \eta)P(5 \nmid \eta) \cdots P(\ell(p) \nmid \eta)}{P(2 \nmid n) P(3 \nmid n)P(5 \nmid n) \cdots P(\ell(p) \nmid n)}. \]
A second heuristic step is then to estimate $P(\ell \nmid \eta)$ by
\[ 1 - \frac{1}{\ell - 1} \text{ if $\ell \nmid p-1$}, \quad \text{and} \quad 1 - \frac{\ell}{\ell^2 - 1} \text{ if $\ell \mid p-1$} \]
(following Theorem~\ref{Lenstra}), and $P(\ell \nmid n)$ by
\[ 1 - \frac{1}{\ell}.\]
One finds
\[ c_p' = \frac{\prod_{\ell \nmid p-1} \left( 1 - \frac{1}{\ell - 1} \right) \cdot \prod_{\ell \mid p-1} \left(1 - \frac{\ell}{\ell^2 - 1} \right)
}{\prod \left( 1 - \frac{1}{\ell} \right)},\] where the products
are over all primes $\ell \leq \ell(p)$ satisfying the
stated conditions. Rearranging the expression shows that
\[ \lim_{p \rightarrow \infty} \left(c_p - c_p'\right) = 0, \]
where $c_p$ is the factor appearing in Conjecture~\ref{galbraithmckee}.

It is tempting to validate the heuristics using an independence argument based on the Chinese Remainder Theorem (for $n$)
and Howe's generalization of Lenstra's theorem (for $\eta$, see \cite{Howe}).
However, this is too naive. By Mertens' theorem and the Prime Number Theorem
\[ \prod_{\ell \leq \sqrt{p} + 1} \left( 1 - \frac{1}{\ell} \right) \approx \frac{2e^{-\gamma}}{\log p} \approx 2e^{-\gamma}P_1(p).\]
Here, $\gamma \approx 0.57722$ is the Euler-Mascheroni constant ($2e^{-\gamma} \approx 1.12292$).
For the heuristics to be justified, we should hence have that
\[ \prod_{\ell \nmid p-1, \ell \leq \sqrt{p} + 1} \left( 1 - \frac{1}{\ell - 1} \right) \cdot \prod_{\ell \mid p-1, \ell \leq \sqrt{p} + 1} \left(1 - \frac{\ell}{\ell^2 - 1} \right) \approx 2e^{-\gamma}P_2(p).\]
With this in mind, the heuristics becomes very subtle: why would both naive estimates be
equally wrong, as Galbraith and McKee call it? We cannot give a satisfying answer, but note the following.
\emph{(i)} The constant $2e^{-\gamma}$, which reflects the ignored dependency between
being divisible by distinct primes, is accumulated in the tail of the product,
with respect to which $\eta$ and $n$ behave much alike.
Stated alternatively, the `local ratios' $P(\ell \nmid \eta)/P(\ell \nmid n)$ converge quickly to $1$.
By considering $c_p$ as the limiting product of these local ratios, rather than the ratio of two diverging products,
one gets a more comfortable underpinning of the conjectured heuristics.
\emph{(ii)} The heuristics is supported by Galbraith and McKee's first
argument in favor of Conjecture~\ref{galbraithmckee}, which uses
different methods (namely, the analytic Hurwitz-Kronecker class number formula).
\emph{(iii)} As far as computationally feasible, the conjectures that we obtain assuming
this principle are confirmed
by experiment in Section~\ref{experimental}.
\emph{(iv)} The constant from Lemma~\ref{galbraithmckeeaverage} provably appeared
in the dual setting of a fixed elliptic curve over $\mathbb{Q}$ reduced
modulo varying primes $p$, see \cite[Theorem~1]{BCD}.

We end this section with a proof of Lemma~\ref{galbraithmckeeaverage}.\\

\noindent \textsc{Proof of Lemma~\ref{galbraithmckeeaverage}.}
First, let us give a heuristic derivation. Let $\ell$ be a prime number. By Dirichlet's theorem, the proportion of
primes $p$ satisfying $\ell \mid p-1$ is $1/(\ell-1)$.
Averaging
out Lenstra's result then gives
\[ P( \ell \mid \eta ) \approx \frac{1}{\ell - 1} \frac{\ell}{\ell^2- 1} + \frac{\ell - 2}{\ell - 1} \frac{1}{\ell - 1} = \frac{\ell^2 - 2}{(\ell^2 - 1)(\ell - 1)}.\]
So
\[ \frac{P(\ell \nmid \eta)}{P(\ell \nmid n)} \approx 1 - \frac{\ell^2 - \ell - 1}{(\ell^2 - 1)(\ell - 1)^2}, \]
and applying the above heuristics yields the requested formula.

To make the argument precise, pick any $\varepsilon >0$. It is easy to see that
there is a uniform bound $L$ such that $| c_p^L - c_p | < \varepsilon/3$ for all $p$ -- where $c_p^L$
is defined as in Conjecture~\ref{galbraithmckee}, but with the product restricted to
primes $\ell$ that do not exceed $L$ -- and such that, similarly,
\[ \left| \, \prod_{\ell \leq L} \left( 1 - \frac{\ell^2 - \ell - 1}{(\ell^2 - 1)(\ell - 1)^2 } \right) \, - \, \prod_{\ell} \left( 1 - \frac{\ell^2 - \ell - 1}{(\ell^2 - 1)(\ell - 1)^2 } \right) \, \right| < \varepsilon/3. \]
However, by the Dirichlet equidistribution of primes, and because we are taking finite products now,
there is an $N$ such that $n \geq N$ implies
\[ \left| \, \frac{1}{\pi(n)-2} \sum_{3<p\leq n} c_p^L \, - \prod_{\ell \leq L} \left( 1 - \frac{\ell^2 - \ell - 1}{(\ell^2 - 1)(\ell - 1)^2 } \right) \, \right| < \varepsilon/3.\]
Combining the three bounds concludes the proof.
\hfill $\blacksquare$

\section{The random matrix model} \label{randommatrices}

\subsection{The genus $1$ case}

Lenstra's Theorem~\ref{Lenstra} can be understood from the
following random matrix point of view. Let $\mathbb{F}_q$ be a
finite field.
Let $N$ be a positive integer
coprime to $q$, and consider the set
\[ \text{GL}_2^{(q)}(\mathbb{Z}/(N)) = \left\{ \ \left. M \in \text{GL}_2(\mathbb{Z}/(N)) \ \right| \ \det M = q \ \right\}. \]
This set is acted upon by $\text{GL}_2(\mathbb{Z}/(N))$, by
conjugation. To any elliptic curve $E / \mathbb{F}_q$, we can
unambiguously associate an orbit of this action by collecting the
matrices of $q$th power Frobenius, considered as an endomorphism
of the $\mathbb{Z}/(N)$-module $E[N]$ of $N$-torsion points,
with respect to all possible bases. Denote this orbit by $\mathcal{F}_E$.

Take $\text{char} \, \mathbb{F}_q > 3$. For any union of orbits $\mathcal{C} \subset \text{GL}_2^{(q)}(\mathbb{Z}/(N))$, let
$P(\mathcal{F}_E \subset \mathcal{C})$ denote the probability
that the orbit associated to the elliptic curve $y^2 = x^3 + Ax + B$,
where $(A,B) \in \mathbb{F}_q$
is chosen from $\mathcal{H}_{AB}$
uniformly at random, is contained in $\mathcal{C}$.

\begin{principle} \label{randommatrixelliptic}
There exist $C_1 \in \mathbb{R}_{>0}$ and $c \in \mathbb{Z}_{> 0}$,
such that
\[ \left | \, P(\mathcal{F}_E \subset \mathcal{C}) - \frac{\# \mathcal{C}}{\#  \emph{GL}_2^{(q)}(\mathbb{Z}/(N)) } \, \right| \leq C_1 N^c/ \sqrt{q}\]
for all choices of $q$, $N$, and $\mathcal{C}$ as above.
\end{principle}

\noindent We use the word `Principle', because, to
our knowledge, no complete proof of this statement has appeared in the literature.
Nevertheless, it is commonly accepted and extensively confirmed by experiment.
It is generally believed to follow from the work of Katz and Sarnak \cite[Theorem~9.7.13]{KatzSarnak}. A strategy of proof
was communicated to us by Katz, and essentially matches with
the approach of Achter \cite[Theorem~3.1]{achter}, who proved Principle~\ref{randommatrixelliptic}
under certain mild restrictions on $q$ and $N$ (using $c=3$). However,
a more classically flavored proof of Principle~\ref{randommatrixelliptic} can presumably be obtained by applying
Chebotarev's density theorem \cite[Proposition~6.4.8]{Fried}
to the function field extension $\mathbb{F}_q(j) \subset \mathbb{F}_q(\zeta_{N})(j) \subset \mathbb{F}_q(\zeta_{N})(X(N))$, where
$\zeta_N$ is a primitive $N$th root of unity, and the latter extension corresponds to the modular cover $X(N) \rightarrow X(1)$, which is known to be defined
over $\mathbb{F}_q(\zeta_{N})$. This approach is currently being elaborated in \cite{CastryckHubrechts}.

%A more classically-flavored approach, based on Chebotarev's
%density theorem applied to the modular cover $X(N) \rightarrow X(1)$ is currently being elaborated
%in \cite{castryckhubrechts}.

Principle~\ref{randommatrixelliptic} indeed allows one to rediscover the asymptotics of
Theorem~\ref{Lenstra}, by counting the matrices $M \in
\text{GL}_2^{(p)}(\mathbb{F}_\ell)$ satisfying $p + 1 -
\text{Tr}(M) = 0$. We leave this as an exercise.

\subsection{The general case} \label{randommatrixgenus2}

Let $\mathbb{F}_q$ and $N$ be as before, and let
$\overline{\mathbb{F}}_q$ be an algebraic closure of
$\mathbb{F}_q$. Let $C/\mathbb{F}_q$ be a complete nonsingular
curve of genus $g \geq 1$ and denote by $A = \text{Jac}(C)$ its
Jacobian. Then $q$th power Frobenius defines an endomorphism of
the $2g$-dimensional $\mathbb{Z}/(N)$-module $A[N]$ of $N$-torsion
points on $A$. Instead of considering all bases, we can make a
more canonical choice by restricting to symplectic bases.
We briefly review how this works.

We employ
the following notation and terminology. For any $n \in \mathbb{N}$, $\mathbb{I}_n$
denotes the $n \times n$ identity matrix, and
$\Omega$ denotes the $2g \times 2g$ matrix
\[ \begin{pmatrix} 0 & \mathbb{I}_g \\ - \mathbb{I}_g & 0 \\ \end{pmatrix}.\]
The group
\[ \text{Sp}_{2g}(\mathbb{Z}/(N)) = \left\{ \left. \ M \in \text{GL}_{2g}(\mathbb{Z}/(N)) \ \right| \
{}^t M \Omega M = \Omega \ \right\} \]
is called the group of symplectic $2g \times 2g$ matrices, and
\[ \text{GSp}_{2g}(\mathbb{Z}/(N)) = \left\{ \left. \ M \in \text{GL}_{2g}(\mathbb{Z}/(N)) \ \right| \
\exists \, d \in \mathbb{Z}/(N) \text{ such that } {}^t M \Omega M = d \Omega \ \right\} \]
is referred to as the group of symplectic similitudes.
It is naturally partitioned into the sets
\[ \text{GSp}_{2g}^{(d)}(\mathbb{Z}/(N)) = \left\{ \left. \ M \in \text{GL}_{2g}(\mathbb{Z}/(N)) \ \right| \
{}^t M \Omega M = d \Omega \ \right\} \]
with $d$ ranging over $\left( \mathbb{Z}/(N) \right)^\times$. An element of
$\text{GSp}_{2g}^{(d)}(\mathbb{Z}/(N))$ is called $d$-symplectic.
Note that $1$-symplectic and symplectic are synonymous.
A classical trick using the Pfaffian shows that the determinant of a symplectic
matrix is $1$. Hence the determinant of a $d$-symplectic matrix is $d^g$.

Symplectic matrices pop up in the study of skew-symmetric, nondegenerate
bilinear pairings on, in our case, $2g$-dimensional $\left( \mathbb{Z}/(N) \right)$-modules.
Such pairings are often called symplectic forms.
For any choice of basis, one can consider the
standard symplectic form $\langle \cdot, \cdot \rangle$, defined by the rule
\[ \langle v,w \rangle = {}^tv \, \Omega \, w. \]
Given any symplectic form, one can always
choose a basis
with respect to which it becomes the standard symplectic form: such a basis is called a symplectic basis
or a Darboux basis. When switching between two symplectic bases
corresponding to the same symplectic form,
the matrix of base change is symplectic, and conversely.

Now for each primitive $N^\text{th}$ root of unity $\zeta_N \in
\overline{\mathbb{F}}_q$, the Weil pairing
\[ e_N : A[N] \times A[N] \rightarrow \langle \zeta_N \rangle, \]
when composed with the (non-canonical) map
\[ \langle \zeta_N \rangle \rightarrow \mathbb{Z}/(N) : \zeta_N^i \mapsto i,\]
is a skew-symmetric and nondegenerate bilinear pairing on $A[N]$.
A corresponding symplectic basis $P_1, \dots, P_g, Q_1, \dots, Q_g$ is
characterized by the properties
\[ e_N (P_i, Q_j) = \zeta_N^{\delta_{ij}}, \quad e_N(P_i,P_j) = e_N(Q_i,Q_j) = 1\]
for all $i,j \in \{1, \dots, g\}$, where $\delta_{ij}$ is the Kronecker symbol.
Because of the $\text{Gal}(\overline{\mathbb{F}}_q,\mathbb{F}_q)$-invariance of the Weil pairing, one has
that
\[ e_N(P^\sigma, Q^\sigma) = e_N(P,Q)^q \]
where $P,Q$ are arbitrary points of $A[N]$ and $\sigma$ is $q$th power Frobenius. Then bilinearity
implies that the matrix $F$ of $\sigma$ with respect to $P_1, \dots, P_g, Q_1, \dots, Q_g$
satisfies
\[ {}^t F \Omega F = q \Omega,\]
i.e.\ $F$ is $q$-symplectic.

As mentioned above, a different choice of symplectic basis yields a matrix obtained from $F$ by
$\text{Sp}_{2g}(\mathbb{Z}/(N))$-conjugation. Next, if $\zeta_N$ is replaced by
another $N^\text{th}$ root of unity $\zeta_{N}^j$, $j \in \left( \mathbb{Z}/(N) \right)^\times$,
then $P_1, \dots, P_g, [j]Q_1, \dots, [j]Q_g$ is a symplectic basis,
and the matrix of Frobenius is $d_j F d_j^{-1}$, where
\[ d_j = \begin{pmatrix} \mathbb{I}_g & 0 \\ 0 & j \mathbb{I}_g \\ \end{pmatrix}. \]
Since $\text{Sp}_{2g}(\mathbb{Z}/(N))$ and the matrices $d_j$ generate
$\text{GSp}_{2g}(\mathbb{Z}/(N))$, we conclude that
we can unambiguously associate to $C$
an orbit of $\text{GSp}_{2g}^{(q)}(\mathbb{Z}/(N))$
under $\text{GSp}_{2g}(\mathbb{Z}/(N))$-conjugation.

We are now ready to formulate the hyperelliptic curve analogue of
Principle~\ref{randommatrixelliptic}. Let $\text{char} \, \mathbb{F}_q > 2$ and $g \geq 1$.
For any union of $\text{GSp}_{2g}(\mathbb{Z}/(N))$-orbits $\mathcal{C} \subset \text{GSp}_{2g}^{(q)}(\mathbb{Z}/(N))$, let
$P(\mathcal{F}_f \subset \mathcal{C})$ denote the probability
that the orbit associated to the complete nonsingular model of the (hyper)elliptic curve $y^2 = f(x)$,
where $f(x) \in \mathbb{F}_q[x]$
is chosen from $\mathcal{H}_{2g+2}$
uniformly at random, is contained in $\mathcal{C}$.

\begin{principle} \label{randommatrix}
There exist $C_1 \in \mathbb{R}_{>0}$ and $c \in \mathbb{Z}_{> 0}$,
such that
\[ \left | \, P(\mathcal{F}_f \subset \mathcal{C}) - \frac{\# \mathcal{C}}{\#  \emph{GSp}_{2g}^{(q)}(\mathbb{Z}/(N)) } \, \right| \leq C_1 N^c/ \sqrt{q} \]
for all choices of $q$, $N$, and $\mathcal{C}$ as above, provided that
$N$ is odd as soon as $g > 2$.
\end{principle}

The condition $N$ odd is due to the fact that we restrict to hyperelliptic
curves, which as soon as $g > 2$ behave non-randomly with respect to $2$-torsion -- see Section~\ref{lenstrageneralization}. If
instead we considered Jacobians of arbitrary curves (e.g.\ in the sense of Section~\ref{asymptotics}),
we expect that this condition can be dropped.

Again we use the word `Principle', because
no complete proof of this statement has appeared in the literature to date.
But again, this presumably follows from the work of Katz and Sarnak \cite[Theorem~9.7.13]{KatzSarnak},
as elaborated by Achter \cite[Theorem~3.1]{achter} under mild restrictions on $q$ and $N$. In his case,
the exponent reads $c= 2g^2 + g$. Achter's result is sufficiently general for
many of our needs below. In particular, it is sufficient for generalizing Theorem~\ref{Lenstra}
to (hyper)elliptic curves of arbitrary genus $g \geq 1$, which is done in Section~\ref{lenstrageneralization}. Also
note that Achter uses $\mathcal{H}_{2g+2}^\text{m}$ rather than $\mathcal{H}_{2g+2}$.

\section{Counting matrices with eigenvalue $1$} \label{countingmatrices}

For use in Sections~\ref{lenstrageneralization} and~\ref{cyclicity}, we study the following general question: given a prime power $q$,
a prime $\ell \nmid q$, an integer $g \geq 0$, and $d \in \{0, \dots, 2g\}$, what
is the proportion $\mathfrak{P}(q,\ell,g,d)$ of matrices in $\text{GSp}_{2g}^{(q)}(\mathbb{F}_\ell)$
for which the eigenspace for eigenvalue $1$ is $d$-dimensional?
The lemma below transfers this question to the
classical groups $\text{Sp}_{2g}(\mathbb{F}_\ell)$ and
$\text{GL}_g(\mathbb{F}_\ell)$. Let $\mathfrak{P}_\text{Sp}(\ell,g,d)$ be
the proportion of matrices in $\text{Sp}_{2g}(\mathbb{F}_\ell)$ having
a $d$-dimensional eigenspace for eigenvalue $1$, and let
$\mathfrak{P}_\text{GL}(\ell,g,d)$ be the corresponding proportion for the general linear group $\text{GL}_g(\mathbb{F}_\ell)$,
where of course $\mathfrak{P}_\text{GL}(\ell,g,d) = 0$ as soon as $d > g$. We include $g=0$ because of the recursive
nature of the arguments below. In this, we assume
that $\text{GSp}_0^{(q)}(\mathbb{F}_\ell) = \text{Sp}_0(\mathbb{F}_\ell) = \text{GL}_0(\mathbb{F}_\ell)$ contains a unique matrix, and that its $1$-eigenspace is
$0$-dimensional. In particular,
$\mathfrak{P}(q,\ell,0,0) = \mathfrak{P}_\text{Sp}(\ell,0,0) = \mathfrak{P}_\text{GL}(\ell,0,0)$ is understood to be $1$.

\begin{lemma} \label{reducetoclassicalgroups}
If $q \equiv 1 \bmod \ell$, then $\mathfrak{P}(q,\ell,g,d) = \mathfrak{P}_\emph{Sp}(\ell,g,d)$. If $q \not \equiv 1 \bmod \ell$,
then $\mathfrak{P}(q,\ell,g,d) = \mathfrak{P}_\emph{GL}(\ell,g,d)$.
\end{lemma}

\noindent \textsc{Proof.} The first statement is a tautology. So assume that $q \not \equiv 1 \bmod \ell$. We follow
ideas of Achter and Holden \cite[Lemma~3.1]{AchterHolden}, which in turn build upon work of Chavdarov \cite{chavdarov}.

First, for $r=0, \dots, g$, let $S(q,\ell,r,d)$ be the subset of $\text{GSp}_{2r}^{(q)}(\mathbb{F}_\ell)$ consisting of those
matrices having characteristic polynomial $(x-1)^r(x-q)^r$ and whose $1$-eigenspace has dimension $d$. Similarly, let $S_\text{GL}(\ell,r,d)$ be the subset of
$\text{GL}_r(\mathbb{F}_\ell)$ consisting of the matrices having characteristic polynomial $(x-1)^r$ and whose $1$-eigenspace has dimension $d$.

We will prove that
\begin{equation} \label{uniformSpGL}
  \# S(q,\ell,r,d) = \frac{ \# \text{Sp}_{2r}(\mathbb{F}_\ell) }{\# \text{GL}_r(\mathbb{F}_\ell)} \cdot \# S_\text{GL}(\ell,r,d).
\end{equation}
By Jordan-Chevalley decomposition, every element $B \in S(q, \ell, r, d)$ can be uniquely written as the commuting product
of a semisimple matrix $B_s$ and a unipotent matrix $B_u$. Necessarily, $B_s \in \text{GSp}_{2r}^{(q)}(\mathbb{F}_\ell)$
has characteristic polynomial $(x-1)^r(x-q)^r$ and $B_u \in \text{Sp}_{2r}(\mathbb{F}_\ell)$ has
characteristic polynomial $(x-1)^{2r}$. By \cite[Lemma~3.4]{chavdarov}, two such matrices $B_s$
must be conjugated by an element of $\text{Sp}_{2r}(\mathbb{F}_\ell)$.
It follows that for fixed $B_s$, the number of corresponding $B$'s in $S(q,\ell,r,d)$ is always
the same.
Since one instance of $B_s$ is $\text{diag}(1,1,\dots,1,q,q, \dots, q)$, whose centralizer
in $\text{Sp}_{2r}(\mathbb{F}_\ell)$
equals
\[ \left\{ \left. \, \begin{pmatrix} M & 0 \\ 0 & {}^t(M^{-1}) \\ \end{pmatrix} \, \right| \, M \in \text{GL}_r(\mathbb{F}_\ell) \, \right\}, \]
the number of possibilities for $B_s$ is $(\# \text{Sp}_{2r}(\mathbb{F}_\ell))/(\# \text{GL}_r(\mathbb{F}_\ell))$,
and for each $B_s$ there are $S_\text{GL}(\ell,r,d)$ appropriate choices for $B_u$. The claim follows.

Now, let $T(q,\ell,g,d)$ be the set
of matrices of $\text{GSp}_{2g}^{(q)}(\mathbb{F}_\ell)$ having a $d$-dimensional $1$-eigenspace,
thus $\# T(q,\ell,g,d) = \mathfrak{P}(q,\ell,g,d) \cdot \# \text{Sp}_{2g}(\mathbb{F}_\ell)$.
We will count the elements $M \in T(q,\ell,g,d)$ separately for each value of $r$, the order
of vanishing at $1$ of the characteristic polynomial $f_M$ of $M$. To $M$ one can associate
a decomposition of the standard symplectic space $\mathbb{F}_\ell^{2g}, \langle \cdot, \cdot \rangle$
of the form $U_{2r} \oplus V_{2(g-r)}$, where $U_{2r}$ and $V_{2(g-r)}$
are $M$-invariant symplectic subspaces of dimensions $2r$ and $2(g-r)$, respectively, satisfying
$f_{M|_{U_{2r}}} = (x-1)^r(x-q)^r$ and $f_{M|_{V_{2(g-r)}}}(1) \neq 0$. Then
\[ \# T(q,\ell,g,d) \, = \, \sum_{r=0}^g \frac{\# \text{Sp}_{2g}(\mathbb{F}_\ell)}{\# \text{Sp}_{2r}(\mathbb{F}_\ell) \cdot \# \text{Sp}_{2(g-r)}(\mathbb{F}_\ell)} \cdot \# S(q,\ell,r,d) \cdot \# T(q,\ell,g-r,0),\]
where the first factor corresponds to the number of ways of decomposing $\mathbb{F}_\ell^{2g}, \langle \cdot, \cdot \rangle$,
the second factor counts the number of possible actions of $M$ on $U_{2r}$, and the third factor counts
the number of actions of $M$ on $V_{2(g-r)}$. We conclude
\begin{equation} \label{recursionGSpq}
 \mathfrak{P}(q,\ell,g,d) \, = \, \sum_{r=0}^g \frac{\# S(q,\ell,r,d) }{\# \text{Sp}_{2r}(\mathbb{F}_\ell)} \cdot \mathfrak{P}(q,\ell,g-r,0).
\end{equation}
Along with
\begin{equation} \label{sumGSpq}
  \sum_{d=0}^g \mathfrak{P}(q,\ell,g,d) = 1,
\end{equation}
one sees that, given the values $\# S(q,\ell,r,d)$, the recursive equation (\ref{recursionGSpq})
determines all $\mathfrak{P}(q,\ell,g,d)$ by induction on $g$: first one determines
$\mathfrak{P}(q,\ell,g,1), \dots, \mathfrak{P}(q,\ell,g,g)$, during which one should use
that $\# S(q,\ell,0,d) = 0$ as soon as $d>0$, and then one uses (\ref{sumGSpq}) to
obtain $\mathfrak{P}(q,\ell,g,0)$.

The statement then follows by noting that one similarly has
\[
 \mathfrak{P}_\text{GL}(\ell,g,d) \, = \, \sum_{r=0}^g \frac{\# S_\text{GL}(\ell,r,d) }{\# \text{GL}_{r}(\mathbb{F}_\ell)} \cdot \mathfrak{P}_\text{GL}(\ell,g-r,0),
\]
along with the same initial conditions. Thus by (\ref{uniformSpGL}),
the probabilities $\mathfrak{P}(q,\ell,g,d)$ and $\mathfrak{P}_\text{GL}(\ell,g,d)$
are solutions to the same recursive equation. By uniqueness they must coincide. \hfill $\blacksquare$\\

Now for the classical groups $\text{Sp}_{2g}(\mathbb{F}_\ell)$ and $\text{GL}_g(\mathbb{F}_\ell)$, these
proportions have been computed before. Parts of the following result
have been (re)discovered by several people (see e.g. \cite{achter2,CohenLenstra}), but the first to obtain
closed formulas for both $\mathfrak{P}_\text{Sp}(\ell,g,d)$ and $\mathfrak{P}_\text{GL}(\ell,g,d)$
seem to be Rudvalis and Shinoda, in an unpublished work of 1988 \cite{RShi} that was reported upon by Fulman \cite{fulman,fulman2}
and, more recently, Lengler \cite{lengler} and Malle \cite{Malle}.

\begin{theorem} \label{everyoneabit}
One has
\begin{align*}
\mathfrak{P}_\emph{GL}(\ell,g,d) & = \frac{1}{\# \emph{GL}_d(\mathbb{F}_\ell)} \cdot \sum_{j=0}^{g-d} \frac{(-1)^j \ell^\frac{j^2-j}{2}}{\ell^{dj} \cdot \# \emph{GL}_j(\mathbb{F}_\ell)},\\
\lim_{g \rightarrow \infty} \mathfrak{P}_\emph{GL}(\ell,g,d) & = \frac{\ell^{-d^2}}{\prod_{j=1}^d (1 - \ell^{-j})^2} \cdot \prod_{j=1}^\infty \left( 1 - \ell^{-j} \right),\\
\mathfrak{P}_\emph{Sp}(\ell,g,d) & = \frac{1}{\# \emph{Sp}_{2k}(\mathbb{F}_\ell)} \cdot \sum_{j=0}^{g-k} \frac{ (-1)^j \ell^{j^2 + j}}{\ell^{2jk} \cdot \# \emph{Sp}_{2j}(\mathbb{F}_\ell)} \quad \text{if $d=2k$ is even,}\\
\mathfrak{P}_\emph{Sp}(\ell,g,d) & = \frac{1}{\ell^{2k+1} \cdot \# \emph{Sp}_{2k}(\mathbb{F}_\ell)} \sum_{j=0}^{g-k-1} \frac{(-1)^j \ell^{j^2 + j}}{\ell^{2j(k+1)} \cdot \# \emph{Sp}_{2j}(\mathbb{F}_\ell)} \quad \text{if $d=2k+1$ is odd,} \\
\lim_{g \rightarrow \infty} \mathfrak{P}_\emph{Sp}(\ell,g,d) & = \frac{\ell^{-\frac{d(d+1)}{2}}}{\prod_{j=1}^d (1 - \ell^{-j})} \cdot \prod_{j=1}^\infty \left( 1 + \ell^{-j} \right)^{-1}.
\end{align*}
\end{theorem}

\noindent \textsc{Proof.} Proofs can be found in \cite[Theorem~6]{fulman} (for everything on the general linear group), and in \cite[Corollary~1]{fulman2}
(for the closed formulas for $\mathfrak{P}_\text{Sp}(\ell,g,d)$) and \cite[Proposition~3.1]{Malle} (for the limit of the latter).
The proofs of Fulman \cite{fulman,fulman2} use the cycle index method, for which, in the symplectic case,
the author assumes that $\ell$ is odd. However, in the meantime, the required theory on cycle indices has
been extended to arbitrary characteristic \cite{fulman3}. The original proof of Rudvalis and Shinoda \cite{RShi} uses integer partitions
and works in full generality. \hfill $\blacksquare$\\

Along with the well-known identities
\begin{equation} \label{cardinalititiesclassicalgroups}
 \# \text{GL}_{g}(\mathbb{F}_\ell) = \ell^\frac{g^2 - g}{2} \prod_{j=1}^g (\ell^j - 1) \quad \text{and}
\quad \# \text{Sp}_{2g}(\mathbb{F}_\ell) = \ell^{g^2} \prod_{j=1}^g (\ell^{2j} - 1)
\end{equation}
(see e.g.\ \cite[Formula~(2.9) and Theorem~3.2]{kim}), Lemma~\ref{reducetoclassicalgroups} and Theorem~\ref{everyoneabit} yield explicit formulas
for each $\mathfrak{P}(q,\ell,g,d)$.

Since the work of Rudvalis and Shinoda cannot be easily accessed, for sake of self-containedness we include
an independent computation of $\mathfrak{P}(q,\ell,g,d)$ for the case where $d=0$. For the purposes of this article, this is the
most prominent case, as we will see in Section~\ref{lenstrageneralization} below.
At the end of this section, we will study the convergence behavior
for $g \rightarrow \infty$ in additional detail.

Is is convenient to consider instead $\mathfrak{Q}(q,\ell,g) = 1 - \mathfrak{P}(q,\ell,g,0)$, the proportion of matrices of $\text{GSp}_{2g}^{(q)}(\mathbb{F}_\ell)$ for which $1$ does appear as an eigenvalue. We prove:
\begin{theorem}\label{Pg} With notation as above, for $g\geq 0$ we have
\begin{align} \label{Pgclose} \mathfrak{Q}(q,\ell,g) = \begin{cases}
 \displaystyle - \sum_{r=1}^g\ell^r \prod_{j=1}^r (1-\ell^{2j})^{-1} & \textnormal{if \ } \ell \mid q-1, \\
\displaystyle - \sum_{r=1}^g \prod_{j=1}^r (1-\ell^j)^{-1} & \textnormal{if \ }  \ell \nmid q-1.
 \end{cases}\end{align}
\end{theorem}
\noindent \textsc{Proof.} Our starting point is the following recursion formula due to Achter and Holden \cite[Lemma~3.2]{AchterHolden}, the proof of which
was our source of inspiration for Lemma~\ref{reducetoclassicalgroups} above: one has
\[ \mathfrak{Q}(q,\ell,g) = \sum_{r=1}^g \frac{S(q, \ell, r)}{ \# \text{Sp}_{2r}(\mathbb{F}_\ell)}
\left( 1 - \mathfrak{Q}(q,\ell,g-r) \right), \]
where
\[ S(q,\ell,r) = \left\{ \begin{array}{ll} \ell^{2r^2} & \text{ if $\ell \mid q-1$,} \\
\ell^{r^2-r} \frac{\# \text{Sp}_{2r}(\mathbb{F}_\ell)}{\# \text{GL}_r(\mathbb{F}_\ell)} & \text{ if $\ell \nmid q-1$} \\
\end{array} \right.\]
and $\mathfrak{Q}(q,\ell,0) = 0$. Clearly, this determines all $\mathfrak{Q}(q,\ell,g)$ uniquely.
Using (\ref{cardinalititiesclassicalgroups}), this can be rewritten as
\begin{align}\label{Prec}
\mathfrak{Q}(q,\ell,g) = \begin{cases}\displaystyle \sum_{r=1}^g \ell^{r^2} ( 1 - \mathfrak{Q}(q, \ell, g-r)) \prod_{j=1}^r ( \ell^{2j} - 1 )^{-1}  &
\quad \text{ if $\ell \mid q - 1$,}
\\
\displaystyle \sum_{r=1}^g \ell^{(r^2-r)/2}  ( 1 - \mathfrak{Q}(q, \ell, g-r)) \prod_{j=1}^r (\ell^j - 1 )^{-1} &
\quad \text{ if $\ell \nmid q - 1$.} \end{cases}
\end{align}
We will prove by induction on $g$ that (\ref{Pgclose}) indeed solves the recursion. We only consider the case
$\ell \nmid q-1$ (the necessary adaptations for the case $\ell \mid q-1$ are straightforward). Define $P_r:=\prod_{j=1}^r(1-\ell^j)^{-1}$ for $r\geq 0$. After rearranging terms and using the induction hypothesis for $g-1$ one finds with some trivial computations that it suffices to prove
\begin{equation}\label{eqPg}
-P_g=\ell^{\frac{g(g-1)}2}\cdot(-1)^g\cdot P_g+\sum_{r=1}^{g-1}\ell^{\frac{r(r-1)}{2}}\cdot(-1)^r\cdot P_r\cdot P_{g-r}.
\end{equation}
We are left with showing that with
\[S_k:=\sum_{r=0}^kT_r\qquad\text{where }T_r:=(-1)^r\cdot\ell^{\frac{r(r-1)}2}\cdot P_r\cdot P_{g-r},\]
we have $S_g=0$. This however follows from the observation that
\[S_k = (-1)^k\cdot\ell^\frac{k(k+1)}2\cdot P_k\cdot P_{g-k}\cdot (1-\ell^g)^{-1}\cdot (1-\ell^{g-k})\]
which can be shown easily using induction on $k$. Indeed: then $S_g=0$ because its last factor is zero. \hfill $\blacksquare$\\

Next, we study the limiting behavior of
$\mathfrak{Q}(q, \ell, g)$ as $g \rightarrow \infty$. Define
\begin{align*}
 \mathcal{E}(q,\ell,g) := \begin{cases}
\displaystyle1-\prod_{j=1}^\infty\left(1+\frac{1}{\ell^j}\right)^{-1} \!\!\!\!\!- \mathfrak{Q}(q,\ell,g) & \textnormal{if \ } \ell \mid q-1, \\
\displaystyle1-\prod_{j=1}^\infty{\left(1-\frac{1}{\ell^j}\right)} - \mathfrak{Q}(q,\ell,g) & \textnormal{if \ } \ell \nmid q-1.  \end{cases}
\end{align*} Then:
\begin{theorem}\label{Pinf} With notation as above, we have
\begin{align*}
\lim_{g\to \infty} \mathfrak{Q}(q,\ell,g) = \begin{cases}
\displaystyle1-\prod_{j=1}^\infty\left(1+\frac{1}{\ell^j}\right)^{-1} & \textnormal{if \ } \ell \mid q-1, \\
\displaystyle1-\prod_{j=1}^\infty{\left(1-\frac{1}{\ell^j}\right)} & \textnormal{if \ } \ell \nmid q-1.  \end{cases} \end{align*}  Moreover, this convergence is alternating, that is, $$\lim_{g\to \infty} \mathcal{E}(q,\ell,g) = 0, {\textnormal{ \ and \ } } (-1)^g\mathcal{E}(q,\ell,g) > 0$$ for each $g\geq 0$.
\end{theorem}

\noindent \textsc{Proof.}
We make use of the well-known $\q$-identity
\begin{align}\label{GRid}\sum_{n\geq 0} \frac{\q^{\frac{n(n-1)}{2}}x^n}{(\q;\q)_n} = \prod_{k=0}^\infty (1+x\q^k)\end{align}
(see for example \cite[II.2]{GR}). Here $(a; \q)_n := \prod_{j=0}^{n-1}(1 - a\q^j)$ is the Pochhammer symbol. We point
out the distinction between $q$ (whose role is limited to separating the cases $\ell \mid q-1$ and $\ell \nmid q-1$) and the
variable $\q$ used here. It is not hard to show that (\ref{Pgclose}) is equivalent to
\begin{align*}
\mathfrak{Q}(q,\ell,g) = \begin{cases}
\displaystyle-\sum_{r=1}^g{ \frac{ \q^{r^2} (-1)^r}{(\q^2;\q^2)_r} }  & \textnormal{if \ } \ell \mid q-1, \\
\displaystyle-\sum_{r=1}^g{ \frac{ \q^{r(r+1)/2} (-1)^r}{(\q;\q)_r} }  & \textnormal{if \ } \ell \nmid q-1,  \end{cases} \end{align*}
where $\q = \ell^{-1}$.

If $\ell \nmid q-1$, it immediately follows that
\begin{align}\lim_{g\to \infty} \mathfrak{Q}(q,\ell,g) = - \sum_{r=1}^\infty \frac{\q^{\frac{r(r+1)}{2}}(-1)^r}{(\q;\q)_r}
= 1-\prod_{n=1}^\infty (1-\q^n) = 1-\prod_{n=1}^\infty \left(1-\frac{1}{\ell^n}\right),\label{qid}\end{align}
where we used (\ref{GRid}) with $x=-q$. To show the convergence is alternating, we have by definition of $\mathcal{E}(q,\ell,g)$ and Theorem \ref{Pg}, that
\begin{align}\label{Eg}
\mathcal{E}(q,\ell,g)= -\sum_{r=g+1}^\infty \prod_{j=1}^r (1-\ell^j)^{-1},
\end{align} which tends to $0$ as $g\to \infty$.  We observe that consecutive summands in (\ref{Eg}) add to \begin{align} -\frac{(-1)^r}{\prod_{j=1}^r (\ell^j-1)} - \frac{(-1)^{r+1}}{\prod_{j=1}^{r+1} (\ell^j-1)} = \frac{(-1)^{r+1}(\ell^{r+1}-2)}{\prod_{j=1}^{r+1} (\ell^j-1)}. \label{consesum}
\end{align}  Now $r\geq 1$ and $\ell$ is prime so that (\ref{consesum}) is positive if and only if $(-1)^{r+1} >0$, which holds if and only if $r$ is odd.  The sum in (\ref{Eg}) begins with an odd index if and only if $g$ is even or $g=0$, which shows that $(-1)^g\mathcal{E}(q,\ell,g)>0$.

If $\ell \mid q-1$, we conclude similarly that
\begin{align*}\lim_{g\to \infty} \mathfrak{Q}(q,\ell,g) &= - \sum_{r=1}^\infty \frac{\q^{r^2}(-1)^r}{(\q^2;\q^2)_r}
= 1-\prod_{n=1}^\infty (1-\q^{2n-1}) \\ &= 1-\prod_{n=1}^\infty (1+\q^n)^{-1} =
1-\prod_{n=1}^\infty \left(1+\frac{1}{\ell^n}\right)^{-1},\end{align*}
by replacing $\q$ by $\q^2$, and setting $x=-\q$.
To show the convergence is alternating, we have by definition of $\mathcal{E}(q,\ell,g)$ and Theorem \ref{Pg}, that
\begin{align}\label{EgB}
\mathcal{E}(q,\ell,g)= -\sum_{r=g+1}^\infty \ell^r \prod_{j=1}^r (1-\ell^{2j})^{-1},
\end{align} which tends to $0$ as $g\to \infty$.  We observe that consecutive summands in (\ref{EgB}) add to \begin{align} -\frac{(-\ell)^r}{\prod_{j=1}^r (\ell^{2j}-1)} - \frac{(-\ell)^{r+1}}{\prod_{j=1}^{r+1} (\ell^{2j}-1)} = \frac{(-\ell)^{r+1}(\ell^{2r+2}-1-\ell)}{\prod_{j=1}^{r+1} (\ell^{2j}-1)}. \label{consesumB}
\end{align}  Again because $r\geq 1$ and $\ell$ is prime, we find that (\ref{consesumB}) is positive if and only if $(-1)^{r+1} >0$ so by the argument given in the previous case when $\ell \nmid q-1$, we have that $(-1)^g\mathcal{E}(q,\ell,g)>0$ in this case as well.
\hfill $\blacksquare$

\section{A generalization of Lenstra's theorem} \label{lenstrageneralization}

With Principle~\ref{randommatrix} in mind, generalizing Lenstra's Theorem~\ref{Lenstra} boils down
to counting matrices $M \in \text{GSp}_{2g}^{(q)}(\mathbb{F}_\ell)$ having $1$ as an eigenvalue.
Indeed, the Jacobian of a curve $C / \mathbb{F}_q$ will have a rational $\ell$-torsion point
if and only if Frobenius acting on $\text{Jac}(C)[\ell]$ has a fixed point, i.e.\ an eigenvector with eigenvalue $1$.

More formally, for every positive integer $g \geq 1$, and
for each pair of distinct primes $p > 2$ and~$\ell$,
let $P(p,\ell,g)$ be the probability that the Jacobian of the (hyper)elliptic curve
$y^2 = f(x)$, with $f(x) \in \mathbb{F}_p[x]$ uniformly randomly chosen from $\mathcal{H}_{2g+2}$,
has rational $\ell$-torsion. Assume that $\ell$ is odd. Then according to Principle~\ref{randommatrix},
there exist $C_1 \in \mathbb{R}_{>0}$ and $c \in \mathbb{Z}_{>0}$, independent of $p$ and $\ell$ (but depending on $g$), such that
\[ \left| P(p,\ell,g) - \mathfrak{Q}(p,\ell,g) \right| \leq C_1 \ell^c / \sqrt{p},\]
where $\mathfrak{Q}(p,\ell,g)$ is defined as in Section~\ref{countingmatrices} above.
This can be considered a proven statement: Achter's proof \cite[Theorem~3.1]{achter}
covers the case where $\mathbb{F}_q$ is a large prime field.
Therefore, we conclude:
\begin{theorem} \label{lenstrageneralized}
There exist $C_1 \in \mathbb{R}_{>0}$ and $c \in \mathbb{Z}_{>0}$, such that
\begin{align*}
\left| P(p,\ell,g) + \sum_{r=1}^g\ell^r \prod_{j=1}^r (1-\ell^{2j})^{-1} \right| & \leq C_1 \ell^c/\sqrt{p} \quad \text{if $\ell \mid p-1$ \quad \text{ and}}\\
\left| P(p,\ell,g) + \sum_{r=1}^g \prod_{j=1}^r (1-\ell^j)^{-1} \right| & \leq C_1 \ell^c/\sqrt{p} \quad \text{if $\ell \nmid p-1$}
\end{align*}
for all pairs of distinct primes $p, \ell >2$.
\end{theorem}
Note once more that $C_1$ and $c$ do depend on $g$.

Theorem~\ref{lenstrageneralized} is invalid for $\ell = 2$:
as soon as $g > 2$, hyperelliptic curves behave unlike general curves with respect to $2$-torsion. But we can
estimate $P(p,2,g)$ using the following slightly simplified result of Cornelissen \cite[Theorem~1.4]{Cornelissen}:
\begin{theorem}[Cornelissen] \label{Cornelissentheorem}
Let $f(x) \in \mathcal{H}_{2g+2}$. Then the Jacobian of the hyperelliptic curve defined by
$y^2 = f(x)$ does \emph{not} have $\mathbb{F}_p$-rational $2$-torsion if and only if
\begin{itemize}
  \item ($g$ odd) $f(x)$ factors as a product of two irreducible polynomials of odd degree;
  \item ($g$ even) $f(x)$ factors as a product of two irreducible polynomials of odd degree, or
  $f(x)$ is irreducible itself.
\end{itemize}
\end{theorem}
Using that a polynomial of degree $d \geq 1$ over $\mathbb{F}_p$ is irreducible with probability
approximately $1/d$, we obtain the following estimates.
\begin{corollary} \label{cornelissencorollary}
If $g$ is odd, then
\[ P(p,2,g) \rightarrow 1 - \sum_{j=0}^{(g-1)/2} \frac{1}{2j+1} \cdot \frac{1}{2g+2 - (2j+1)} \qquad \text{as $p \rightarrow \infty$}, \]
whereas if $g$ is even, we have
\[ P(p,2,g) \rightarrow 1 - \frac{2g}{(2g+2)^2} - \sum_{j=0}^{g/2} \frac{1}{2j+1} \cdot \frac{1}{2g+2 - (2j+1)} \qquad \text{as $p \rightarrow \infty$}. \]
In particular, we have
\[ \lim_{g,p \rightarrow \infty} P(p,2,g) = 1,\]
hence Theorem~\ref{limithyperelliptic} holds.
\end{corollary}
Note again that for $g \in \{1,2\}$, where the random matrix heuristics are assumed to apply (and in fact provably do
for $\ell=2$ --- see Corollary~\ref{randommatrixlevel2} for $g=2$, exercise for $g=1$), we obtain $P(p,2,1) = 2/3$ and $P(p,2,2) \approx 25/46$,
which is the same as if we would have evaluated the second formula of Theorem~\ref{Pg} in $\ell=2$.

We are now ready to derive Conjectures~\ref{ourconjecture} and~\ref{limitconjecture}, and to prove Lemma~\ref{ourconjectureaverage}.\\

\noindent \textsc{Derivation of Conjecture~\ref{ourconjecture}.}
Let $\mathbb{F}_p$ be a large prime field and let $\ell$ be a prime different from its
characteristic $p$. From Theorem~\ref{Pg}, we see that the
probability that the Jacobian of $y^2 = f(x)$, with $f(x)$ chosen from $\mathcal{H}_6$
uniformly at random, has rational $\ell$-torsion is approximately
\[ \frac{ \ell(\ell^4 - \ell - 1) }{(\ell^4-1)(\ell^2-1)} \text{ if $\ell \mid p-1$} \qquad \text{and} \qquad \frac{\ell^2 - 2}{(\ell^2-1)(\ell-1)} \text{ if $\ell \nmid p-1$}.\]
Note that because $g=2$, these limiting probabilities are also valid for $\ell=2$.
Applying the heuristics from Section~\ref{framework} then yields the requested formula for $c_p$.
One new point of concern is that $\ell(p)$, which should now be the largest prime
for which $\ell(p) \leq (\sqrt{p}+1)^2$, exceeds $p$. Therefore,
we should take into account the contribution of $\ell = p$. But since
we take $p \rightarrow \infty$, it suffices that the probability of not having $p$-torsion
tends to $1$. This follows from Principle~\ref{levelp} (Section~\ref{extensionfields}) below.
\hfill $\blacksquare$\\

\noindent \textsc{Proof of Lemma~\ref{ourconjectureaverage}.}
This is entirely analogous to the proof of Lemma~\ref{galbraithmckeeaverage}.\hfill $\blacksquare$\\

\noindent \textsc{Derivation of Conjecture~\ref{limitconjecture}.}
Applying our heuristics, using the probabilities
given in Theorem~\ref{Pinf}, we obtain
\[c_p = \prod_{\ell \nmid p-1} \frac{\prod_{j=1}^\infty \left(1 - \frac{1}{\ell^j} \right)}{1 - \frac{1}{\ell}}
\cdot \prod_{\ell \mid p - 1} \frac{\prod_{j=1}^\infty \left(1 + \frac{1}{\ell^j} \right)^{-1}}{1 - \frac{1}{\ell}} \]
Note that we also use these probabilities for $\ell=2$, since we expect the random matrix statement
from Pinciple~\ref{randommatrix} to apply in arbitrary level $N$ (in the current, more general framework of selecting curves from $\mathcal{M}_g$
uniformly at random).
Rearranging factors gives
\[ c_p = \prod_\ell \prod_{j=2}^\infty \left(1 - \frac{1}{\ell^j} \right) \cdot
\prod_{\ell \mid p - 1} \prod_{j=1}^\infty \left(1 - \frac{1}{\ell^{2j}} \right)^{-1}, \]
from which the requested formula follows.
\hfill $\blacksquare$\\

\noindent We remark that the average setups (Lemmata~\ref{galbraithmckeeaverage} and~\ref{ourconjectureaverage})
can be thought of as taking matrices at random from $\text{GSp}_{2g}(\mathbb{F}_\ell)$, rather
than $\text{GSp}_{2g}^{(p)}(\mathbb{F}_\ell)$.\\

It is interesting to note,
using Theorem~\ref{Pinf}, that as the genus $g$ grows, the average value $\cpave$
oscillates, but converges rapidly to its limiting value. This is illustrated numerically in Table~\ref{T_alternating}.
Of all genera, elliptic curves disfavor prime orders to the biggest extent, and the Jacobians of genus 2 curves disfavor prime orders to the least extent.

\begin{table}
\begin{center}
\begin{tabular}{|l|c|}
\hline
$g$&$\cpave$\\\hline
$1$&0.50516617\\
$2$&0.69463828\\
$3$&0.68851794\\
$4$&0.68857163\\
$5$&0.68857149\\
$6$&0.68857149\\
$7$&0.68857149\\
\hline
\end{tabular}
\caption{Value of $\cpave$ for growing genus, i.e.\ the constants appearing in
Lemmata~\ref{galbraithmckeeaverage} and~\ref{ourconjectureaverage},
and their higher genus analogues.}
\label{T_alternating}
\end{center}
\end{table}

\section{The case of a rational Weierstrass point} \label{rationalweierstrass}

In many applications, often cryptographic, one restricts to genus
$2$ curves of the form $y^2 = f(x)$ where $f(x)$ is chosen from
\[ \mathcal{H}_5^\text{m} = \left\{ f \in \mathbb{F}_q[x] \, | \, \text{$f$ monic and square-free}, \, \deg f = 5 \right\} \]
uniformly at random.
Stated
more geometrically, one restricts to genus $2$ curves having a
rational Weierstrass point. However, the latter description is not free of
ambiguities. Namely, consider the notion of randomness in
which $f(x)$ is taken from
\[ \mathcal{H}_6^{(>0)} = \left\{ f \in \mathbb{F}_q[x] \, | \, \text{$f$ square-free}, \, \deg f = 6, \, \exists \, a \in \mathbb{F}_q : f(a) = 0 \right\} \]
uniformly at random.
Then this is fundamentally different from the $\mathcal{H}_5^\text{m}$-setting. To illustrate this: the
probability that the Jacobian of a randomly chosen curve
has even order tends to $4/5
= 0.8$ with respect to $\mathcal{H}_5^\text{m}$, whereas it tends to $311/455 \approx
0.68$ with respect to $\mathcal{H}_6^{(>0)}$. Both statements will be proven
below.

The main conclusion of this section will be, however, that the distribution
of Frobenius acting on any \emph{odd}-torsion subgroup of the Jacobian
is barely affected by this ambiguity. In Section~\ref{oddlevel},
we will show:
\begin{theorem}~\label{maintheoremdeg5}
Let $N$ be an odd positive integer, let $q$ be an odd prime power coprime to $N$, and let $\mathcal{H}$ be either $\mathcal{H}_5^\emph{m}$, $\mathcal{H}_6^{(>0)}$
or $\mathcal{H}_6$. For any subset $\mathcal{C} \subset \emph{GSp}_4^{(q)}(\mathbb{Z}/(N))$
that is closed under $\emph{GSp}_4(\mathbb{Z}/(N))$-conjugation,
let $P(\mathcal{F}_f \subset \mathcal{C})$ be defined as in Section~\ref{randommatrixgenus2}, where
now $f$ is chosen from $\mathcal{H}$ uniformly at random. If Principle~\ref{randommatrix} holds, then there
exist $C_1 \in \mathbb{R}_{>0}$ and $c \in \mathbb{Z}_{> 0}$ such that
\[ \left|  P \left( \mathcal{F}_f \subset \mathcal{C} \right) - \frac{\# \mathcal{C}}{\# \emph{GSp}_4^{(q)}(\mathbb{Z}/(N))} \right| \leq C_1 N^c / \sqrt{q} \]
for all choices of $q$ and $\mathcal{C}$ as above.
\end{theorem}
For $\mathcal{H}_5^\text{m}$, we remark that it is presumably possible to prove Theorem~\ref{maintheoremdeg5} directly from Katz-Sarnak \cite[Theorem~9.7.13]{KatzSarnak}, i.e.\ independently of Principle~\ref{randommatrix}, in the same way as a proof of Principle~\ref{randommatrix} is expected to work, using that the family corresponding to $\mathcal{H}_5^\text{m}$ has the largest possible monodromy group \cite[10.1.18]{KatzSarnak}.

As an immediate application, one obtains:\\

\noindent \textsc{Heuristic derivation of Conjecture~\ref{imposingweierstrass}.} By Theorem~\ref{maintheoremdeg5},
we only need to replace the factor $\frac{38}{45}$, corresponding to the prime $\ell = 2$, by $\frac{2}{5}$. So the
correcting factor is $\frac{9}{19}$.  \hfill $\blacksquare$

\subsection{Rational $2$-torsion in genus $2$} \label{even}

Some material in this section has appeared in the literature before, see
e.g.\ \cite[Section 2]{vandergeer}.

\begin{lemma} \label{weierstrassgenerates2torsion}
Every non-trivial $2$-torsion point on the Jacobian of a genus $2$ curve
over $\mathbb{F}_q$ (thought
of as a divisor class) contains a unique pair of divisors $\{P_i - P_j, P_j - P_i \}$, where
$P_i$ and $P_j$ are distinct Weierstrass points.
\end{lemma}
\textsc{Proof.}
It is obvious that $P_i - P_j$ and $P_j - P_i$ are linearly equivalent, and that
they map to a $2$-torsion point on the Jacobian. By Riemann-Roch, this point is non-trivial
and two different pairs give rise to distinct $2$-torsion points. Since there are $15$ non-trivial
$2$-torsion points on the Jacobian of a genus 2 curve, and since there are $15$ pairs in a set of $6$ elements,
the correspondence must be $1$-to-$1$.
\hfill $\blacksquare$\\

\noindent We immediately obtain (compare with Theorem~\ref{Cornelissentheorem}):

\begin{lemma} \label{deg5reducible}
The Jacobian of a genus $2$ curve over $\mathbb{F}_q$ defined by an equation of the
form $y^2 = f(x)$ with $f \in \mathcal{H}_5^\emph{m}$ resp.\ $f \in \mathcal{H}_6$ has a
non-trivial rational $2$-torsion point if and only if $f$ is
reducible resp. $f$ has a factor of degree $2$.
\end{lemma}
\textsc{Proof.}
By Lemma~\ref{weierstrassgenerates2torsion}, there exists a non-trivial rational $2$-torsion
point if and only if there are Weierstrass points $P_1$ and $P_2$ such that
$\{P_1,P_2\}$ is closed under $q$th power Frobenius. \hfill $\blacksquare$\\

\noindent This allows us to estimate the probability that the Jacobian has even order.

\begin{lemma} \label{mod2counts}
Let $f_5^\emph{m} \in \mathcal{H}_5^\emph{m}$,
$f_6^{(>0)} \in \mathcal{H}_6^{(>0)}$ and $f_6 \in \mathcal{H}_6$ be
chosen uniformly at random. Let $C_5^\emph{m}$, $C_6^{(>0)}$ and $C_6$ denote the corresponding
genus $2$ curves. Then as $q \rightarrow \infty$
\begin{enumerate}
  \item[(i)] $P( \# \emph{Jac}(C_5^\emph{m})(\mathbb{F}_q) \text{ is even}) \rightarrow 4/5$;
  \item[(ii)] $P( \# \emph{Jac}(C_6)(\mathbb{F}_q) \text{ is even}) \rightarrow 26/45$;
  \item[(iii)] $P( \# \emph{Jac}(C_6^{(>0)})(\mathbb{F}_q) \text{ is even}) \rightarrow 311/455$.
%  \item[(iv)] $P( f_6 \in \mathcal{H}_6^{(w)} ) \rightarrow 455/720$ $(=91/144)$.
\end{enumerate}
\end{lemma}
\textsc{Proof.}
We leave this as an exercise, or refer to Table~\ref{T_2torsion} below. \hfill $\blacksquare$\\

We will now describe the symplectic structure of the $2$-torsion subgroup in more
detail. Fix a genus $2$ curve $C / \mathbb{F}_q$ and let $P_1, \dots, P_6$ be
its Weierstrass points. Following Lemma~\ref{weierstrassgenerates2torsion},
every non-trivial element of $\text{Jac}(C)[2]$ can be identified with a unique pair
of distinct points $\{P_i, P_j\}$, and the group structure can be described by the rules
\[ \left\{ \begin{array}{l} \{P_i,P_j\} + \{P_i,P_j \} = 0 \\ \{P_i,P_j\} + \{P_i,P_k\} = \{P_j,P_k\} \quad
\text{if $j \neq k$} \\ \{P_i,P_j\} + \{ P_k,P_\ell \} = \{ \text{remaining two points} \} \quad \text{if $\{i,j\} \cap \{k,\ell\} = \emptyset$}. \\ \end{array} \right. \]
The Weil pairing can be seen to satisfy
\[ e_2(\{P_i,P_j\}, \{P_k,P_\ell\}) = (-1)^{\# \{i,j,k,\ell\}} \]
for all $i,j,k,\ell \in \{1, \dots, 6\}$.

%Note that in the above deduction, $\text{Sp}_4(\mathbb{Z}/(2))$ was seen as the group of symplectic \emph{transformations}.
%By fixing a symplectic basis, it becomes isomorphic to the group of symplectic matrices. By allowing all symplectic
%bases, a symplectic transformation actually corresponds to a \emph{conjugacy class} of $\text{Sp}_4(\mathbb{Z}/(2))$,
%and hence of $\text{Sym}(6)$.

%In particular, $q$th power Frobenius corresponds to a conjugacy class of $\text{Sym}(6)$, which
%we will call the \emph{Frobenius conjugacy class}. Note that the number of fixed points
%of a permutation is an invariant of its conjugacy class. It is then easy
%to see that for each $r=1, \dots, 6$, the Frobenius conjugacy class has
%exactly $r$ fixed points if and only if $C$ has $r$ $\mathbb{F}_q$-rational
%Weierstrass points.

We use this to prove the following.

\begin{theorem} \label{lemmaW}
Let $q$ be an odd prime power. There exist $\mathcal{W}_0, \dots, \mathcal{W}_6 \subset \emph{Sp}_4(\mathbb{F}_2)$
such that for any curve $C / \mathbb{F}_q$ of genus $2$, any symplectic basis
of $\emph{Jac}(C)[2]$, and any $r \in \{0, \dots, 6\}$,
the matrix $F$ of $q$th power Frobenius with respect to this basis
satisfies
\[ F \in \mathcal{W}_r \quad \text{if and only if $C$ has $r$ rational Weierstrass points}.\]
The cardinalities of the $\mathcal{W}_r$ are $265$, $264$, $135$, $40$, $15$, $0$ and $1$, respectively.
\end{theorem}
\noindent \textsc{Proof.} There exist $6$ subsets $U \subset \text{Jac}(C)[2]$ that are maximal
with respect to the condition that $u_1,u_2 \in U$ and $u_1\neq u_2$ implies $e_2(u_1,u_2) = -1$,
namely
\[ U_i = \left\{ \left. \{P_i,P_j\} \, \right| \, j \in \{1,2, \dots, 6\} \setminus \{i\} \right\} \quad \text{for $i=1, \dots, 6$}. \]
Since $N=2$, the choice of a primitive $N$th root of unity is canonical, hence
the Weil pairing  defines unambiguously a symplectic pairing
on $\text{Jac}(C)[2]$. After having fixed a symplectic basis, every
symplectic matrix induces a permutation of $\{ U_1, \dots, U_6\}$. In fact,
this induces a group isomorphism $\text{Sp}_4(\mathbb{F}_2) \rightarrow \text{Sym}(6)$.
Indeed, it is easy to see that the above induces an injective group homomorphism, and
surjectivity follows from $\#\text{Sp}_4(\mathbb{F}_2) = \#\text{Sym}(6) = 720$.
Then the sets $\mathcal{W}_r$ are
the preimages under this isomorphism
of the set of permutations having
exactly $r$ fixed points. While the isomorphism
depends on the choice of symplectic basis, the sets $\mathcal{W}_r$ do not,
because they are invariant under conjugation.\hfill $\blacksquare$\\

Pushing the argument a little further, one actually
sees that the conjugacy class of Frobenius, which
under the above group isomorphism corresponds to a conjugacy class of $\text{Sym}(6)$,
is completely determined by the factorization pattern of $f(x)$, and conversely. Note that there are $11$ conjugacy classes
in $\text{Sym}(6) \cong \text{Sp}_4(\mathbb{F}_2)$, and that there are $11$ ways to partition the number $6$.
Since the probability of having a certain factorization pattern is easily estimated
using the well-known fact
that a polynomial of degree $d$ is irreducible with probability about $1/d$,
this unveils the complete stochastic picture of $\text{Jac}(C)[2]$, as shown
in Table~\ref{T_2torsion}.
\begin{table}
\begin{center}
\begin{tabular}{|lr|lr|lrrrr|}
\hline
\multicolumn{2}{|c|}{$\mathcal{H}_6$} & \multicolumn{2}{c|}{$\mathcal{H}_5^\text{m}$} & \multicolumn{5}{c|}{conjugacy classes of $\text{Sp}_4(\mathbb{F}_2)$}\\\hline
pattern&prob.&pattern&prob.&representant&size&order&$\mathbb{F}_q$-rank&trace\\\hline & & & & & & & &\vspace{-13pt}\\
6 & $\approx \frac{1}{6}$ &  &  & \tiny $\left(\begin{smallmatrix} 1 & 1 & 0 & 1 \\ 0 & 1 & 0 & 1 \\ 1 & 1 & 1 & 1 \\ 0 & 0 & 1 & 1 \\ \end{smallmatrix}\right)$ & 120 & 6 & 0 & 0 \\[4pt]
5,1 & $\approx \frac{1}{5}$ & 5  & $\approx \frac{1}{5}$ & \tiny $\left(\begin{smallmatrix} 0 & 1 & 1 & 0 \\ 0 & 0 & 1 & 0 \\ 0 & 1 & 0 & 1 \\ 1 & 1 & 0 & 1 \\ \end{smallmatrix}\right)$ & 144 & 5 & 0 & 1 \\[4pt]
4,2 & $\approx \frac{1}{8}$ &  &  & \tiny $\left(\begin{smallmatrix} 0 & 1 & 0 & 0 \\ 1 & 0 & 0 & 0 \\ 0 & 0 & 0 & 1 \\ 1 & 0 & 1 & 0 \\ \end{smallmatrix}\right)$ & 90 & 4 & 1 & 0 \\[4pt]
4,1,1 & $\approx \frac{1}{8}$ & 4,1  & $\approx \frac{1}{4}$ & \tiny $\left(\begin{smallmatrix} 1 & 0 & 1 & 1 \\ 1 & 0 & 0 & 1 \\ 0 & 1 & 0 & 1 \\ 1 & 1 & 1 & 1 \\ \end{smallmatrix}\right)$ & 90 & 4 & 1 & 0 \\[4pt]
3,3 & $\approx \frac{1}{18}$ &  &  & \tiny $\left(\begin{smallmatrix} 1 & 0 & 1 & 1 \\ 0 & 1 & 1 & 0 \\ 0 & 1 & 0 & 0 \\ 1 & 1 & 0 & 0 \\ \end{smallmatrix}\right)$ & 40 & 3 & 0 & 0 \\[4pt]
3,2,1 & $\approx \frac{1}{6}$ & 3,2  & $\approx \frac{1}{6}$ & \tiny $\left(\begin{smallmatrix} 1 & 1 & 0 & 0 \\ 0 & 1 & 1 & 1 \\ 0 & 1 & 0 & 1 \\ 1 & 1 & 1 & 1 \\ \end{smallmatrix}\right)$ & 120 & 6 & 1 & 1 \\[4pt]
3,1,1,1 & $\approx \frac{1}{18}$ & 3,1,1 & $\approx \frac{1}{6}$ & \tiny $\left(\begin{smallmatrix} 1 & 0 & 0 & 0 \\ 0 & 0 & 0 & 1 \\ 0 & 0 & 1 & 0 \\ 0 & 1 & 0 & 1 \\ \end{smallmatrix}\right)$ & 40 & 3 & 2 & 1 \\[4pt]
2,2,2 & $\approx \frac{1}{48}$ &  &  & \tiny $\left(\begin{smallmatrix} 0 & 0 & 0 & 1 \\ 1 & 0 & 1 & 0 \\ 0 & 1 & 0 & 1 \\ 1 & 0 & 0 & 0 \\ \end{smallmatrix}\right)$ & 15 & 2 & 2 & 0 \\[4pt]
2,2,1,1 & $\approx \frac{1}{16}$ & 2,2,1 & $\approx \frac{1}{8}$ & \tiny $\left(\begin{smallmatrix} 0 & 0 & 1 & 0 \\ 0 & 1 & 0 & 0 \\ 1 & 0 & 0 & 0 \\ 0 & 1 & 0 & 1 \\ \end{smallmatrix}\right)$ & 45 & 2 & 2 & 0 \\[4pt]
2,1,1,1,1 & $\approx \frac{1}{48}$ & 2,1,1,1 & $\approx \frac{1}{12}$ & \tiny $\left(\begin{smallmatrix} 0 & 1 & 1 & 0 \\ 0 & 1 & 0 & 0 \\ 1 & 1 & 0 & 0 \\ 1 & 1 & 1 & 1 \\ \end{smallmatrix}\right)$ & 15 & 2 & 3 & 0 \\[4pt]
1,1,1,1,1,1 & $\approx \frac{1}{720}$ & 1,1,1,1,1 & $\approx \frac{1}{120}$ & \tiny $\left(\begin{smallmatrix} 1 & 0 & 0 & 0 \\ 0 & 1 & 0 & 0 \\ 0 & 0 & 1 & 0 \\ 0 & 0 & 0 & 1 \\ \end{smallmatrix}\right)$ & 1 & 1 & 4 & 0 \\[4pt]
\hline
\end{tabular}
\caption{Factorization patterns of $f(x) \in \mathcal{H}_6, \mathcal{H}_5^\text{m}$ and the corresponding Frobenius conjugacy classes.
For instance, the pattern 3,1,1,1 means that $f(x) \in \mathcal{H}_6$ factors into three linear polynomials and one irreducible cubic polynomial.
The probability of this event is approximately $\frac{1}{3} \cdot \frac{1}{3!} \frac{1}{1} = \frac{1}{18}$. The corresponding
conjugacy class of Frobenius is generated by the depicted matrix and contains $40$ elements. Every such
element has order $3$ and trace $1$, and its eigenspace for eigenvalue $1$ is $2$-dimensional (i.e.\ $\dim \text{Jac}(C)[2](\mathbb{F}_q) = 2$).}
\label{T_2torsion}
\end{center}
\end{table}
\begin{corollary} \label{randommatrixlevel2}
Principle~\ref{randommatrix} holds for $g=N=2$.
\end{corollary}
\noindent \textsc{Proof.} This can be read off from the above table. The only additional concern is the bound on the error term, but this is easily verified. \hfill $\blacksquare$

\subsection{Equidistribution in odd level} \label{oddlevel}

In this section, we will prove Theorem~\ref{maintheoremdeg5}.
Consider $f \in \mathcal{H}_6^{(>0)}$, so that $y^2 = f(x)$ defines
a genus $2$ curve having a rational Weierstrass point $(a,0)$.
Then the birational change of variables
\[x \leftarrow \frac{1}{x} + a, \qquad y \leftarrow \frac{y}{x^3}, \]
transforms this into $y^2 = f'(x)$
with $f' \in \mathcal{H}_5$.
This leads us to defining a relation
\[ \rho \subset \mathcal{H}_6^{(>0)} \times \mathcal{H}_5 \]
associating to $f \in \mathcal{H}_6^{(>0)}$ all polynomials of $\mathcal{H}_5$ that
can be obtained through the above procedure.
However, this correspondence is not uniform, because of
the number of choices that can be made for $a$, i.e. the
number of rational roots of $f$. This is the reason why the
notions of randomness with respect to $\mathcal{H}_5$ (or $\mathcal{H}_5^\text{m}$) and
$\mathcal{H}_6^{(>0)}$
are fundamentally different, as reflected in Lemma~\ref{mod2counts}.

We are led to introducing the following notation. For $r \in \{0, \dots, 6\}$, define
\[ \mathcal{H}_6^{(r)} = \left\{ f \in \mathbb{F}_q[x] \, | \, \text{$f$ square-free}, \, \deg f = 6, \, \text{$f$ has precisely $r$ rational zeroes} \right\} \]
so that
\begin{equation} \label{partition6}
 \mathcal{H}_6 = \bigsqcup_{r=0}^6 \mathcal{H}_6^{(r)} \quad \text{and} \quad \mathcal{H}_6^{(>0)} = \bigsqcup_{r=1}^6 \mathcal{H}_6^{(r)}.
\end{equation}
Similarly, for $r \in \{0, \dots, 5\}$ we introduce
\[ \mathcal{H}_5^{(r)} = \left\{ f \in \mathbb{F}_q[x] \, | \, \text{$f$ square-free of degree $5$, $f$ has precisely $r$ rational zeroes} \right\},\]
so that
\[ \mathcal{H}_5 = \bigsqcup_{r=0}^5 \mathcal{H}_5^{(r)}. \]
Note that $\mathcal{H}_6^{(5)}$ and $\mathcal{H}_5^{(4)}$ are empty. We
implicitly omit these sets to avoid probabilities of the type $\frac{0}{0}$.
Similarly, we assume that $q > 6$ so that none of the other sets are empty.

Now because of (\ref{partition6}), to prove Theorem~\ref{maintheoremdeg5} for $\mathcal{H}_6^{(>0)}$, it suffices
to do so for each $\mathcal{H}_6^{(r)}$ ($r=1,\dots, 6$). Similarly, by the discussion in Section~\ref{randomness}
we can use $\mathcal{H}_5$ instead of $\mathcal{H}_5^\text{m}$, and it is sufficient
to prove Theorem~\ref{maintheoremdeg5} for $\mathcal{H}_5^{(r)}$ ($r=0, \dots, 5$) in this case.
Finally, by the lemma below, the cases $\mathcal{H}_5^{(r)}$
can in turn be reduced to the cases $\mathcal{H}_6^{(r)}$.

\begin{lemma} \label{locallyuniform}
Let $S_0 = \left\{ f \in \mathcal{H}_5 \, | \, f(0) \neq 0 \right\}$.
For each $r = 1, \dots, 6$, the restriction of $\rho$ to
\[ \mathcal{H}_6^{(r)} \times \left( \mathcal{H}_5^{(r-1)} \cap S_0 \right) \]
is uniform.
\end{lemma}
\textsc{Proof.} This is immediate. \hfill $\blacksquare$\\

We are now ready to prove Theorem~\ref{maintheoremdeg5}.\\

\noindent \textsc{Proof of Theorem~\ref{maintheoremdeg5}.} By the above discussion, it
suffices to estimate the conditional probabilities
\[ P(\mathcal{F}_f \subset \mathcal{C} \, | \, f \in \mathcal{H}_6^{(r)}) =
\frac{P(\mathcal{F}_f \subset \mathcal{C} \text{ and }
f \in \mathcal{H}_6^{(r)})}{P(f \in \mathcal{H}_6^{(r)})} \]
for $r=1, \dots, 6$.
By Theorem~\ref{lemmaW},
$f \in \mathcal{H}_6^{(r)}$ is equivalent to
saying that the conjugacy class of Frobenius,
acting on the $2$-torsion points of the Jacobian of $y^2 = f(x)$,
is contained in $\mathcal{W}_r$. Denote this conjugacy class
by $\mathcal{F}_{f,2}$. Similarly, let
$\mathcal{F}_{f,2N}$ denote the conjugacy class of Frobenius acting on the $2N$-torsion points.

Since $N$ is odd, we have a canonical isomorphism
\[ \text{GSp}_4^{(q)}(\mathbb{Z}/(2N)) \cong \text{GSp}_4^{(q)}(\mathbb{F}_2) \oplus \text{GSp}_4^{(q)}(\mathbb{Z}/(N)), \]
allowing us to consider $\mathcal{W}_r \oplus \mathcal{C}$ as
a subset of $\text{GSp}_4^{(q)}(\mathbb{Z}/(2N))$. Because
it is the union of a number of orbits under $\text{GSp}_4(\mathbb{Z}/(2N))$-conjugation,
there exist $C_1 \in \mathbb{R}_{>0}$ and $c \in \mathbb{Z}_{>0}$, such that
\begin{equation} \label{level2N}
 \left | \, P(\mathcal{F}_{f,2N} \subset \mathcal{W}_r \oplus \mathcal{C}) - \frac{\# \left(\mathcal{W}_r \oplus \mathcal{C} \right)}{\#  \text{GSp}_4^{(q)}(\mathbb{Z}/(2N)) } \, \right| \leq C_1 N^c/ \sqrt{q}
\end{equation}
for all choices of $q$, $N$ and $\mathcal{C}$.
In particular, for $N=1$ this gives
\begin{equation} \label{level2}
 \left | \, P(f \in \mathcal{H}_6^{(r)}) - \frac{\# \mathcal{W}_r }{\#  \text{GSp}_4^{(q)}(\mathbb{F}_2) } \, \right| \leq C_1 / \sqrt{q}.
\end{equation}
Since
\[ P(\mathcal{F}_{f,2N} \subset \mathcal{W}_r \oplus \mathcal{C}) = P(\mathcal{F}_f \subset \mathcal{C} \text{ and } \mathcal{F}_{f,2} \subset \mathcal{W}_r) = P(\mathcal{F}_f \subset \mathcal{C} \text{ and }
f \in \mathcal{H}_6^{(r)}) \]
and
\[ \frac{\# (\mathcal{W}_r \oplus \mathcal{C})}{\# \text{GSp}_4^{(q)}(\mathbb{Z}/(2N))} =
\frac{\# \mathcal{W}_r}{\# \text{GSp}_4^{(q)}(\mathbb{F}_2)} \cdot \frac{\# \mathcal{C}}{\# \text{GSp}_4^{(q)}(\mathbb{Z}/(N))}, \]
inequality (\ref{level2N}) can be rewritten as
\[ \left | \, P(\mathcal{F}_f \subset \mathcal{C} \, | \, f \in \mathcal{H}_6^{(r)}) -
\frac{ \frac{\# \mathcal{W}_r}{\# \text{GSp}_4^{(q)}(\mathbb{F}_2)} }{P(f \in \mathcal{H}_6^{(r)})}
\cdot \frac{\# \mathcal{C}}{\# \text{GSp}_4^{(q)}(\mathbb{Z}/(N))}
\, \right| \leq \frac{C_1 N^c/ \sqrt{q}}{P(f \in \mathcal{H}_6^{(r)})}.\]
It follows from (\ref{level2}) that
there is a $C_2 \in \mathbb{R}^+$ such that
\[ \left | \, P(\mathcal{F}_f \subset \mathcal{C} \, | \, f \in \mathcal{H}_6^{(r)}) -
 \frac{\# \mathcal{C}}{\# \text{GSp}_4^{(q)}(\mathbb{Z}/(N))}
\, \right| \leq C_2 N^c/ \sqrt{q}\]
for all choices of $q$, $N$ and $\mathcal{C}$. This ends the proof. \hfill $\blacksquare$

\section{The number of points on the curve itself} \label{curveitself}

Up to now we have focused entirely on the number of rational points on the Jacobian of a curve. However, the random matrix framework allows us to consider the number of rational points on the curve itself as well.

For any pair of distinct primes $p>2$ and  $\ell$, and any $t \in \mathbb{F}_\ell$, we define the following constants:
\begin{eqnarray*}
a_{\ell,t,p} &:=& \#\{(x,y)\in\mathbb{F}_\ell^\times\times(\mathbb{F}_\ell^\times\backslash \{-p\}) \ |\ (x+y/x)(1+p/y)=t\},\\
A_{\ell,t,p} &:=& \ell^4((\ell-1)(\ell-2)+a_{\ell,t,p})+\begin{cases}\ell^6-\ell^4 & \textnormal{if $t= 0$,}\\ 0  & \textnormal{otherwise,}\end{cases}\\
B_{\ell} &:=& \ell^4(\ell^2-1)^2,\\
C_{\ell,t} &:=& \ell^5(\ell-1)(\ell^3-\ell-1)+\begin{cases}\ell^7-\ell^6 & \textnormal{if $t= 0$,}\\ 0  & \textnormal{otherwise.}\end{cases}
\end{eqnarray*}
Note that it is probably impossible to find a simple formula for $a_{\ell,t,p}$ since, in general, it describes the number of points on an elliptic curve over $\mathbb{F}_\ell$ (though it is clear that $a_{\ell,t,p}$ lies close to $\ell$).
Let $P(p,\ell,t)$ be the probability that the number of rational points on the nonsingular complete model of
the curve $C : y^2=f(x)$, with $f(x)$ chosen uniformly at random from $\mathcal{H}_6$,
is congruent to $p+1-t$ modulo $\ell$.
\begin{theorem}\label{propositionOnSizeCurve}
There exist $C_1\in\mathbb{R}_{>0}$ and $c \in \mathbb{Z}_{>0}$, such that
\[\left| P(p,\ell,t) - \cfrac{A_{\ell,t,p}+B_{\ell}+C_{\ell,t}}{\ell^4\cdot(\ell^4-1)\cdot(\ell^2-1)} \right|\leq C_1\ell^c/\sqrt{p}\]
for all $p,\ell,t$ as above.
\end{theorem}
\textsc{Proof.} Because the trace of a matrix is invariant under conjugation, it suffices by Principle~\ref{randommatrix} (proven for $\ell$ odd by Achter \cite[Theorem~3.1]{achter}, and for $\ell=2$ in Corollary~\ref{randommatrixlevel2})  to count the number of matrices $M$ in $\text{GSp}_{4}^{(p)}(\mathbb{F}_\ell)$ with trace $t$, and show that it equals $A_{\ell,t,p}+B_{\ell}+C_{\ell,t}$.
Our main tool is the
following Bruhat decomposition of
$\text{Sp}_{4}(\mathbb{F}_\ell)$, proven by Kim \cite{kim}.
Consider the group
\begin{equation}\label{equationBruhat}P = \left\{ \left. \begin{pmatrix} A & AB \\ 0 & {}^t{A}{^{-1}} \\ \end{pmatrix} \, \right|
\, A, B \in \mathbb{F}_\ell^{2 \times 2}, \text{$A$ invertible,
$B$ symmetric} \right\},\end{equation}
then we have the disjoint union
\[ \text{Sp}_4(\mathbb{F}_\ell) = P \ \sqcup \  P
\sigma_1 P
\ \sqcup \ P \sigma_2 P \]
where
\[ \sigma_1 = \begin{pmatrix} 0 & 0 & 1 & 0 \\ 0 & 1 & 0 & 0 \\ -1 & 0 & 0 & 0 \\ 0 & 0 & 0 & 1 \\ \end{pmatrix} \quad
\text{and} \quad \sigma_2 = \Omega = \begin{pmatrix} 0 & \mathbb{I}_2 \\ - \mathbb{I}_2 & 0 \\ \end{pmatrix}. \]
For $r \in \{1, 2\}$, consider the subgroup
\[ A_r = \{ \, M \in P \, | \, \sigma_r M \sigma_r^{-1} \in P \, \}. \]
Then one can find unique representatives for the elements of $P\sigma_rP$ by rewriting
\[  P \sigma_r P = P \sigma_r (A_r \backslash P ), \]
where $A_r \backslash P$ should be seen as a set of representatives
of the right cosets of $A_r$ in $P$. This implies that
\[ | P \sigma_r P | = |P| \cdot |A_r \backslash P|. \]
One can prove (see \cite{kim}) that $|A_1 \backslash P| = \ell^2 + \ell$ and $|A_2 \backslash P| = \ell^3$.
Taking $\sigma_0=\mathbb{I}_4$, the Bruhat decomposition of $\text{Sp}_{4}(\mathbb{F}_\ell)$ implies
the following partition of $\text{GSp}_4^{(p)}(\mathbb{F}_\ell)$:
\[ \text{GSp}_{4}^{(p)}(\mathbb{F}_\ell) = \bigsqcup_{r=0}^2 d_pP \sigma_r P. \]
We will do a component-wise count of the number of matrices having trace $t$.
First we observe that
\[ \left| \left\{ \left. M \in d_pP \sigma_r P \, \right| \, \text{Tr}(M)=t \right\} \right| =
|A_r \backslash P | \cdot \left| \left\{ \left. M \in d_pP \sigma_r \, \right| \, \text{Tr}(M)=t \right\} \right| \]
for $r = 1,2$. Indeed,
every element of $d_pP \sigma_r P$
has a unique representation of the form
\[ d_p M \sigma_r N \]
with $M \in P$ and $N \in A_r \backslash P$ (where $A_r \backslash P$
is thought of as a set of representatives of the right cosets of $A_r$). Using this
representation, the map
\[ d_p P \sigma_r P \rightarrow d_p P \sigma_r : d_p M \sigma_r N \mapsto d_p (d_p^{-1} N d_p M) \sigma_r \]
is surjective and $|A_r \backslash P|$-to-$1$. Since
$d_p M \sigma_r N$ and $d_p (d_p^{-1} N d_p M) \sigma_r$ are conjugated, the observation follows.

A matrix $M\in d_pP$ can be written as $\left(\begin{smallmatrix} A & AB \\ 0 & p\cdot{}^t{A}{^{-1}}\end{smallmatrix}\right)$ with
$A \in \text{GL}_2(\mathbb{F}_\ell)$ and $B \in \mathbb{F}_\ell^{2 \times 2}$ symmetric.

First, we consider $M\sigma_1$, whose trace equals $-(AB)_{1,1}+A_{2,2}+(p\cdot{}^t{A}{^{-1}})_{2,2}$, where the index notation refers to the corresponding entries. Fix $A$ and let $B$ vary. Then because $(AB)_{1,1}=A_{1,1}B_{1,1}+A_{1,2}B_{2,1}$ and not both $A_{1,1}$ and $A_{1,2}$ can be zero, we find that each trace occurs equally often. We conclude
that traces are uniformly distributed in $d_pP\sigma_1$.
Next, for $M\sigma_2$ we find that $\text{Tr}(M\sigma_2)=-\text{Tr}(AB)$, which is uniformly distributed for all $A$ not of the form $\left(\begin{smallmatrix} 0&a\\-a&0\end{smallmatrix}\right)$, and which is zero if $A$ does have this form. Using the above formulas for $| A_r \backslash P |$ and using $|\text{GL}_2(\mathbb{F}_\ell)|=\ell(\ell^2-1)(\ell-1)$,
we find that the number of matrices in $d_pP\sigma_1 \sqcup d_pP\sigma_2$ having trace $t$ equals $B_{\ell}+C_{\ell,t}$.

Finally we consider $M\in d_pP$ when $\text{Tr}(M)=\text{Tr}(A)+\text{Tr}(pA^{-1})$. We write $A=\left(\begin{smallmatrix}a & b\\c & d\end{smallmatrix}\right)$ and let $\delta=ad-bc$ be its determinant. Clearly $\text{Tr}(M)=\text{Tr}(A)\cdot(1+p/\delta)$. There are
$\ell(\ell^2-1)$ matrices $A$ with determinant $-p$, in which case this trace equals 0. So suppose that $\delta\neq -p$. When $a=0$ it is easy to see that we have uniform distribution, so we also suppose that $a\neq 0$. We can replace $d$ by $(\delta+bc)/a$ and again, if $b\neq 0$ we will find uniformity. Finally the case $b=0$ gives as trace
\[(a+\delta/a)(1+p/\delta),\] so that an easy calculation shows that the number of matrices in $d_p P$ with trace $t$ equals $A_{\ell,t,p}$.\hfill$\blacksquare$\\

Table~\ref{T_curveitself} gives the respective probabilities for various small $\ell$.
\begin{table}
\begin{center}
\begin{tabular}{|l|r|r|r|r|r|r|}\cline{1-4}
\vphantom{$\sum^{a^b}$} & $p\bmod \ell\ \backslash\ t$ & 0 & 1 & \multicolumn{3}{|c}{}\\\cline{2-4}
\vphantom{$\sum^{a^b}$}\raisebox{1.5ex}[0cm][0cm]{$\ell=2$} & 1 & $\frac{26}{45}$ & $\frac{19}{45}$ & \multicolumn{3}{|c}{}\\\cline{1-5}

\vphantom{$\sum^{a^b}$} & $p\bmod \ell\ \backslash\ t$ & 0 & 1 & 2 & \multicolumn{2}{|c}{}\\\cline{2-5}
\vphantom{$\sum^{a^b}$}{$\ell=3$}  & 1 & $\frac{46}{128}$ & $\frac{41}{128}$ & $\frac{41}{128}$ & \multicolumn{2}{|c}{}\\
\vphantom{$\sum^{a^b}$} & 2 & $\frac{58}{160}$ & $\frac{51}{160}$ & $\frac{51}{160}$ & \multicolumn{2}{|c}{}\\\cline{1-7}

\vphantom{$\sum^{a^b}$} & $p\bmod \ell\ \backslash\ t$ & 0 & 1 & 2 & 3 & 4\\\cline{2-7}
\vphantom{$\sum^{a^b}$}  & 1 & $\frac{3094}{14976}$ & $\frac{2969}{14976}$ & $\frac{2972}{14976}$ & $\frac{2972}{14976}$ & $\frac{2969}{14976}$\\
\vphantom{$\sum^{a^b}$}{$\ell=5$}  & 2 & $\frac{774}{3744}$ & $\frac{743}{3744}$ & $\frac{742}{3744}$ & $\frac{742}{3744}$ & $\frac{743}{3744}$\\
\vphantom{$\sum^{a^b}$}\vphantom{$\sum^{a^b}$}  & 3 & $\frac{774}{3744}$ & $\frac{742}{3744}$ & $\frac{743}{3744}$ & $\frac{743}{3744}$ & $\frac{742}{3744}$\\
\vphantom{$\sum^{a^b}$} & 4 & $\frac{3094}{14976}$ & $\frac{2972}{14976}$ & $\frac{2969}{14976}$ & $\frac{2969}{14976}$ & $\frac{2972}{14976}$\\\hline
\end{tabular}
\caption{Distribution of Frobenius traces modulo small $\ell$ for $y^2 = f(x)$, with
$f(x) \in \mathcal{H}_6$ chosen at random.}
\label{T_curveitself}
\end{center}
\end{table}
Note that the probabilities of $C$ resp.\ $\text{Jac}(C)$ having an even number of rational points
are the same, despite the fact that these events do not coincide.
Also note from Table~\ref{T_curveitself} that trace $0$ is favored. This is a general phenomenon
that can be seen as follows.
It is not hard to verify that if $2t(t^2-16p)\equiv 0\bmod \ell$, the curve
$(x + y/x)(1 + p/y) = t$ in the definition of
$a_{\ell,t,p}$ is reducible or has genus $0$,
in which case $a_{\ell,t,p}$ can be explicitly computed. It is equal to zero if $\ell=2$. For $t\equiv 0\bmod\ell$ and $\ell>2$ we can compute the following estimate for $P(p,\ell,t)$:
\[
\frac{\ell^9-\ell^6-\ell^5-\ell^4}{\ell^4(\ell^4-1)(\ell^2-1)} = \frac{\ell^3-\ell-1}{(\ell^2-1)^2}
\]
if $p$ is a square modulo $\ell$ and $$\frac{\ell^9-\ell^6-\ell^5+\ell^4}{\ell^4(\ell^4-1)(\ell^2-1)} = \frac{\ell^3+\ell-1}{\ell^4-1}$$ otherwise. Both probabilities are indeed larger than $1/\ell$.
If $p\equiv t^2/16\bmod \ell$ and hence $t\not\equiv 0\bmod \ell$ we obtain
$$\frac{\ell^9-\ell^7-\ell^6-\ell^5-\ell^4}{\ell^4(\ell^4-1)(\ell^2-1)} = \frac{\ell^5-\ell^3-\ell^2-\ell-1}{(\ell^4-1)(\ell^2-1)} $$ if $\ell\equiv 1\bmod 4$ and finally when $\ell\equiv 3\bmod 4$ we find $$\frac{\ell^9-\ell^7-\ell^6-\ell^5+\ell^4}{\ell^4(\ell^4-1)(\ell^2-1)} = \frac{\ell^5-\ell^3-\ell^2-\ell+1}{(\ell^4-1)(\ell^2-1)}.$$\\

\noindent \textsc{Heuristic derivation of Conjecture~\ref{curveitselfconj}.} The number
of rational points on the curve defined by $y^2 = f(x)$ is divisible by $\ell$ if and only if
its trace $t$ is congruent to $p+1$ mod $\ell$. Thus, by Theorem~\ref{propositionOnSizeCurve}, the probability that this number of points
is not divisible by $\ell$ can be estimated by
\[ \frac{\beta_{\ell,p}}{(\ell^4 - 1)(\ell^2-1)}, \]
where $\beta_{\ell,p}$ is as in the introductory Section~\ref{curveitselfintro}. Dividing by $1 - \frac{1}{\ell}$ and
taking the product then gives the constant $c_p$ from Conjecture~\ref{curveitselfconj}. The factor corresponding to
$\ell = 2$ can be read off from the table above (or from Table~\ref{T_2torsion}).
When switching from $\mathcal{H}_6$ to $\mathcal{H}_5$, following Theorem~\ref{maintheoremdeg5} and using Table~\ref{T_2torsion},
we should replace the factor $\frac{38}{45}$ by $\frac{16}{15}$. \hfill $\blacksquare$

\section{The probability of cyclicity} \label{cyclicity}

In this section, we will estimate the probability $P(p,g)$ that the group of rational points of the Jacobian of the (hyper)elliptic curve
$C:y^2 = f(x)$, with $f(x)$ chosen from
$\mathcal{H}_{2g+2}$
uniformly at random, is cyclic. This question is of a different type from what we have considered so far.
We use the following heuristic reasoning. Note that $\text{Jac}(C)(\mathbb{F}_p)$ is cyclic if and only
if $\text{Jac}(C)[\ell](\mathbb{F}_p)$ is cyclic for each prime $\ell$. The probabilities of the latter events
can be estimated using Principle~\ref{randommatrix}: for each $\ell \neq p$, this is approximately
\[ \mathfrak{P}(p,\ell,g,0) + \mathfrak{P}(p,\ell,g,1),\]
where the notation from Section~\ref{countingmatrices} is used. For a reason similar to the one explained
in the derivation of Conjecture~\ref{ourconjecture} in Section~\ref{lenstrageneralization}, we will omit the contribution of $\ell = p$.
Then the idea is to assume independence
and naively multiply these proportions. As suggested by our experiments in Section~\ref{experimental},
this gives very accurate predictions for $g \in \{1,2\}$. In particular, an effect of the type reflected in
Mertens' theorem seems absent in this non-relative setting.
For $g=1$, the heuristics confirm a formula proven by Vl\u{a}du\c{t} \cite[Theorem~6.1]{vladut}.\\

\noindent \textsc{Heuristic derivation of Conjecture~\ref{cyclicitygenus2}.} The formulas of Theorem~\ref{everyoneabit} for $g=2$ give
\[ \mathfrak{P}(p,\ell,2,0) + \mathfrak{P}(p,\ell,2,1) = \left\{ \begin{array}{lr}
\frac{\ell^8 - \ell^6 - \ell^5 - \ell^4 + \ell^2 + \ell + 1}{\ell^2(\ell^4-1)(\ell^2-1)} & \text{if $\ell \mid p-1$,} \\
1 - \frac{1}{\ell(\ell^2-1)(\ell-1)} & \text{if $\ell \nmid p-1$.} \\  \end{array} \right.\]
Multiplying gives the conjectured formula. If we switch from $\mathcal{H}_6$ to $\mathcal{H}_5^\text{m}$,
the leading factor $\frac{151}{180}$ should be replaced by $\frac{37}{60}$, as can be read off from
Table~\ref{T_2torsion}.\hfill $\blacksquare$\\

\noindent \textsc{Proof of Theorem~\ref{limitcyclichyperelliptic}.} This is analogous to the proof of Theorem~\ref{limithyperelliptic}
(see Corollary~\ref{cornelissencorollary}). In fact, the original version of Cornelissen's Theorem~\ref{Cornelissentheorem} \cite[Theorem~1.4]{Cornelissen}
is much stronger and describes the rank of $\text{Jac}(C)[2](\mathbb{F}_p)$
in terms of the factorization pattern of $f(x)$. E.g., it suffices
that $f(x)$ has at least $4$ distinct factors for the rank to be at least $2$. From this, one verifies that for $g \rightarrow \infty$,
this rank will be $2$ or larger with a probability converging to $1$. \hfill $\blacksquare$\\

\noindent \textsc{Heuristic derivation of Conjecture~\ref{limitcyclic}.} This is a combination of the derivations
of Conjectures~\ref{limitconjecture} and~\ref{cyclicitygenus2}, the details of which we leave to the reader. \hfill $\blacksquare$\\

As in the case of primality, we list the average values (in the sense of Conjecture~\ref{averagingoverp})
of the probabilities of cyclicity for growing genus in Table~\ref{T_cyclicity}. Again one notices that the convergence is alternating (although
we did not elaborate the details of a proof of this) and fast.
\begin{table}
\begin{center}
\begin{tabular}{|l|r|}
\hline
$g$&factor\\\hline
$1$&0.81375191\\
$2$&0.80882586\\
$3$&0.80924272\\
$4$&0.80923674\\
$5$&0.80923677\\
$6$&0.80923677\\
$7$&0.80923677\\
\hline
\end{tabular}
\caption{Average conjectured probability of being cyclic for growing genus.}
\label{T_cyclicity}
\end{center}
\end{table}

\section{Extension fields} \label{extensionfields}

In this section, we briefly discuss how our heuristics can be adapted
to the setting of finite fields $\mathbb{F}_{p^k}$ of growing extension degree, over a fixed prime field $\mathbb{F}_p$.
In this situation one can no longer neglect the contribution of the prime $\ell = p$.

Let $C/ \mathbb{F}_{p^k}$ be a complete nonsingular curve of genus $g \geq 1$ and, as before, denote by $A = \text{Jac}(C)$
its Jacobian. One has
\[ A[p] \cong \left(\mathbb{F}_p \right)^r \]
for some $0 \leq r \leq g$.
We assume that if $k$ is large and one picks $C$ at random (e.g.\ from
\[ \mathcal{M}_g = \left\{ \, \text{curves of genus $g$ over $\mathbb{F}_{p^k}$} \, \right\} / \cong_{\mathbb{F}_{p^k}} \]
uniformly at random), one has $r=g$ with probability $\approx 1$. This is reasonable, because
the moduli space $\mathcal{A}_g$ of abelian varieties of dimension $g$ is stratified by rank, the
stratum corresponding to $r=g$ having the biggest dimension \cite[Theorem~4.1]{NormanOort}.
We do not claim a proof of this assumption however, although for hyperelliptic curves
this is a known fact \cite{achterpries,prieszhu}.
If $r=g$, then the matrix
of $p^k$th power Frobenius acting on $A[p]$ with respect to any
$\mathbb{F}_p$-basis is an element of $\text{GL}_g(\mathbb{F}_p)$.
Thus, in that case, we can unambiguously associate to $C$ a conjugacy class of
matrices of $p^k$th power Frobenius, denoted by $\mathcal{F}_C$.
The expectation is that for every union of conjugacy classes $\mathcal{C} \subset \text{GL}_g(\mathbb{F}_p)$,
the probability that $\mathcal{F}_C \subset \mathcal{C}$ becomes proportional to $\# \mathcal{C}$ (as $k \rightarrow \infty$).

Returning to hyperelliptic curves, let $P(\mathcal{F}_{f,h} \subset \mathcal{C})$
be the probability that the conjugacy class of Frobenius associated to the hyperelliptic curve
$y^2 + h(x)y = f(x)$, where $(f,h)$ is chosen
from $\mathcal{H}_{g+1,2g+2}$
uniformly at random, is contained in $\mathcal{C}$. As explained in Section~\ref{randomness}, for $p > 2$ one can
assume $h(x) = 0$ and $f(x)$ chosen from $\mathcal{H}_{2g+2}$ if wanted.
\begin{principle} \label{levelp}
Let $g \in \{1,2\}$. There exist $C_1 \in \mathbb{R}_{>0}$ and $c \in \mathbb{Z}_{>0}$ such that
\[ \left| P(\mathcal{F}_{f,h} \subset \mathcal{C}) - \frac{\# \mathcal{C}}{\# \emph{GL}_g(\mathbb{F}_p)} \right| \leq C_1 p^c / \sqrt{p^k} \]
for all choices of $p, k$, and $\mathcal{C}$ as above.
\end{principle}
The assumption $g \in \{1,2\}$ is a `safety' measure, because we do not feel comfortable with the behavior of the hyperelliptic locus inside
$\mathcal{A}_g$ as soon as $g > 2$. In fact, even for $g=2$ some prudence is needed with respect to Principle~\ref{levelp}:
the literature seems to contain much less evidence in its favor than in the cases of Principle~\ref{randommatrixelliptic} and Principle~\ref{randommatrix}.

In contrast, for $g=1$ Principle~\ref{levelp}
can be proven by applying the Hasse-Weil bound to the Igusa curve $\text{Ig}(p)$,
whose $\mathbb{F}_{p^k}$-rational points essentially parameterize pairs $(E,P)$, where $E / \mathbb{F}_{p^k}$ is an elliptic
curve and $P \in E[p](\mathbb{F}_{p^k})$. A more elementary but longer proof is given
below. We include it because we believe some intermediate statements are interesting in their own right
(in fact, we develop a version of \cite[Theorem~V.4.1]{Silverman}, which is on the Legendre family, for Weierstrass equations).
First note that Principle~\ref{levelp} is trivial for $p=2$ and for $p=3$,
in the latter case because quadratic twisting provides
a bijection between the set of elliptic curves having trace $1$
mod $3$ and the set of elliptic curves with trace $2$ mod $3$.

\begin{theorem} \label{theoremmodp}
Let $p \geq 5$ be a prime number, let $k \geq 1$ be an integer,
and let $t \in \{1, \dots, p-1\}$. Let $S_t$ be the set of couples
in
\[ S = \mathcal{H}_{A,B} = \left\{ \left. (A,B) \in (\mathbb{F}_{p^k})^2 \, \right| \, 4A^3 + 27B^2 \neq 0 \right\} \]
for which the trace $T$ of the $p^k$th power Frobenius of the elliptic curve given by
$y^2 = x^3 + Ax + B$
satisfies $T \equiv t \mod p$.
Then $\#S = p^{2k} - p^k$ and
\[ \left| \# S_t - \frac{\#S}{p-1} \right| \leq 3p^{\frac{3k}{2} + 1}.\]
\end{theorem}
\textsc{Proof.} We leave it as an exercise to show that $\#S = p^{2k} - p^k$.

For each $(A,B) \in S$, one has that $T$ mod $p$ equals the norm
(with respect to $\mathbb{F}_{p^k} / \mathbb{F}_p$) of the
coefficient $c_{A,B}$ of $x^{p-1}$ in
\[ \left(x^3 + Ax + B\right)^{\frac{p-1}{2}} \]
(see the proof of \cite[Theorem V.4.1(a)]{Silverman}).
Lemma~\ref{squarefree} below shows that for every $\gamma \in
\mathbb{F}_{p^k}^\times$, the polynomial $c_{A,B} -
\gamma$ is absolutely irreducible when $A$ and $B$ are considered to be variables.

Now write $S'_t$ for the set of couples $(A,B) \in
(\mathbb{F}_{p^k})^2$ in which $c_{A,B}$ evaluates to an element
$\gamma \in \mathbb{F}_{p^k} \setminus \{0\}$ with norm $t$
(regardless of the condition $4A^3 + 27B^2 \neq 0$).
There are
\[ \frac{p^k - 1}{p-1} \]
such $\gamma$'s. For each of these the polynomial $c_{A,B} -
\gamma$ defines a plane affine curve, by the claimed irreducibility. Its degree is bounded by $d = 3(p-1)/2$, hence its
(geometric) genus is at most $(d-1)(d-2)/2$, and the number of
points at infinity is at most $d$. Therefore the set $S'_\gamma
\subset S'_t$ of couples satisfying $c_{A,B} = \gamma$ is subject
to
\[ \left| \# S'_\gamma - (p^k + 1) \right| \leq (d-1)(d-2)\sqrt{p^k} + d \leq \frac{9}{4}p^{\frac{k}{2} + 2}\]
by the Hasse-Weil bound. Note that $c_{A,B} = \gamma$ defines
an affine, possibly singular curve, so some caution is needed when applying
the Hasse-Weil bound. See \cite[Theorem~5.4.1]{Fried} for the details.

Summing up, and using $(p^k-1)/(p-1) \leq \frac{5}{4}p^{k-1}$ (since $p \geq 5$),
\[ \left| \# S'_t - \frac{p^{2k} - 1}{p-1} \right| \leq \frac{45}{16}p^{\frac{3k}{2} + 1}.\]
Because $\#(S'_t \setminus S_t) \leq p^k$ and $5p^{k-1} \leq p^k \leq \frac{1}{11}p^{\frac{3}{2}k + 1}$,
we obtain
\[ \left| \# S_t - \frac{p^{2k} - p^k}{p-1} \right| \leq \left| \# S_t - \frac{p^{2k} - 1}{p-1} \right| + \frac{p^k - 1}{p-1}
\leq \left(\frac{45}{16} + \frac{1}{11} + \frac{5}{4}\cdot \frac{1}{55} \right) p^{\frac{3k}{2} + 1},\]
which ends the proof. \hfill $\blacksquare$

\begin{lemma} \label{squarefree}
Let $p\geq 5$ be a prime number and let $c_{A,B} \in \mathbb{F}_p[A,B]$ be the
coefficient of $x^{p-1}$ in
\[ \left(x^3 + Ax + B \right)^{\frac{p-1}{2}} \quad \in \mathbb{F}_p[A,B][x].  \]
Then $c_{A,B}$ is homogeneous of $(2,3)$-weighted degree
$(p-1)/2$, nonzero, and absolutely squarefree. As a consequence,
for any $\gamma \in \overline{\mathbb{F}}_p^\times$, the polynomial
\[ c_{A,B} - \gamma \quad \in \overline{\mathbb{F}}_p[A,B]\] is
irreducible.
\end{lemma}
\noindent \textsc{Proof.} One verifies that
\begin{equation} \label{sum}
 c_{A,B} = \sum_{ i = \left\lceil \frac{p-1}{6} \right\rceil}^{\left\lfloor \frac{p-1}{4} \right\rfloor} { \frac{p-1}{2} \choose i}
{i \choose {3i - \frac{p-1}{2}}} A^{3i - \frac{p-1}{2}} B^{\frac{p-1}{2} - 2i}
\end{equation}
from which it immediately follows that $c_{A,B}$ is nonzero and homogeneous of degree $(p-1)/2$ if
we equip $A$ and $B$ with weights $2$ and $3$ respectively. It is easy to verify that $A$ and $B$ appear as a factor at most once.

Let $c_{A,B}'$ be obtained from $c_{A,B}$ by deleting the factors $A$ and $B$ when possible. Define
$\varepsilon_A$ (resp. $\varepsilon_B$) to be $1$ if a factor $A$ (resp. $B$) was deleted, and $0$ otherwise.
Then $c_{A,B}'$ is still homogeneous, of degree $(p-1)/2 - 2\varepsilon_A - 3\varepsilon_B$.
After dividing by a suitable power of $A$ and considering
the resulting polynomial in the single variable $B^2/A^3$, one verifies that $c'_{A,B}$ splits (over $\overline{\mathbb{F}}_p$)
\begin{equation} \label{factorizationcab}
c(B^2 - a_1A^3)(B^2 - a_2A^3) \cdots (B^2 - a_rA^3)
\end{equation}
with $r = \frac{1}{6}\left((p-1)/2 - 2\varepsilon_A - 3\varepsilon_B \right)$ and all $c,a_i \neq 0$.
Each of these factors corresponds to a $j_i \neq 0,1728$ for which
the elliptic curve over $\overline{\mathbb{F}}_p$ with $j$-invariant $j_i$ is supersingular,
and conversely all supersingular $j$-invariants different from $0$ and $1728$
must be represented this way. Now the number of supersingular $j$-invariants
different from $0$ and $1728$ is precisely given by $r$
(see the proof of \cite[Theorem V.4.1(c)]{Silverman}). Therefore,
all factors in (\ref{factorizationcab}) must be different, and in particular $c_{A,B}$ must be
squarefree.

Now let $\gamma \in \overline{\mathbb{F}}_p^\times$ and suppose we had a
nontrivial factorization
\[ c_{A,B} - \gamma = (F_1 + X_1)(F_2 + X_2), \]
where $F_1$ and $F_2$ are the components of highest (weighted) degree of the
respective factors. Then it follows that $F_1F_2 = c_{A,B}$, so
$F_1$ and $F_2$ cannot have a common factor. It also follows that
\begin{equation}\label{vanish}
X_1F_2 + X_2F_1 + X_1X_2 + \gamma = 0.
\end{equation}
Let $X_1'$ and $X_2'$ be the components of highest degree of $X_1$
and $X_2$ respectively. Suppose $\deg X_1F_2 > \deg X_2F_1$. Then $X_1'F_2$ is zero, because it cannot be cancelled in
(\ref{vanish}). But then $X_1' = X_1 = 0$ and we run into a
contradiction. By symmetry, we conclude that $\deg X_1F_2 = \deg
X_2F_1$. But then $X_1'F_2 + X_2'F_1 = 0$. So all factors of $F_1$
must divide $X_1'F_2$, which is impossible unless $X_1' = 0$, and
we again run into a contradiction. \hfill $\blacksquare$\\

We end this section with a derivation of Conjectures~\ref{ext1} and~\ref{ext2}.
To apply our heuristics, we need to generalize the material from Section~\ref{lenstrageneralization}.
In analogy with the notation employed there, for any prime power $q$, any
prime number $\ell$, and any integer $g \geq 1$, let $P(q,\ell,g)$ be the probability that
the Jacobian of the hyperelliptic curve $y^2 + h(x)y = f(x)$, with $(f,h)$ chosen from
$\mathcal{H}_{g+1,2g+2}$
uniformly at random, has an $\mathbb{F}_q$-rational $\ell$-torsion point.
Let $\mathfrak{Q}(q,\ell,g)$ be the proportion of matrices of $\text{GSp}_{2g}^{(q)}(\mathbb{F}_\ell)$
having $1$ as an eigenvalue if $\ell \nmid q$, and the proportion of
matrices of $\text{GL}_g(\mathbb{F}_\ell)$ if $\ell \mid q$. Then according to Principles~\ref{randommatrix}
and~\ref{levelp}, if $g\in \{1,2\}$ we have that $P(q,\ell,g) \rightarrow \mathfrak{Q}(q,\ell,g)$ as $q \rightarrow \infty$.
Recall from Theorem~\ref{Pg} that one has
\begin{align} \label{Pggeneralclose} \mathfrak{Q}(q,\ell,g) = \begin{cases}
 \displaystyle - \sum_{r=1}^g\ell^r \prod_{j=1}^r (1-\ell^{2j})^{-1} & \textnormal{if \ } \ell \mid q-1, \\
\displaystyle - \sum_{r=1}^g \prod_{j=1}^r (1-\ell^j)^{-1} & \textnormal{if \ }  \ell \nmid q-1
 \end{cases}\end{align}
if $\ell \nmid q$. However, the same formula applies for $\ell \mid q$, because in case $\ell \nmid q-1$,
the proportion of matrices of $\text{GSp}^{(q)}_{2g}(\mathbb{F}_\ell)$ having $1$ as an eigenvalue
equals the corresponding proportion for $\text{GL}_g(\mathbb{F}_\ell)$ anyway, due to Lemma~\ref{reducetoclassicalgroups}.
In other words, one can blindly adapt Theorem~\ref{Pg} to this more general setting.
Therefore:\\

\noindent \textsc{Heuristic derivation of
Conjectures~\ref{ext1} and~\ref{ext2}.} This is a copy of the heuristic derivations
of Conjectures~\ref{galbraithmckee} and~\ref{ourconjecture}. \hfill $\blacksquare$

\section{Experimental evidence} \label{experimental}

The following tables present experimental data in support of Conjectures~\ref{galbraithmckee}--\ref{cyclicitygenus2}.
Table~\ref{T_galbraithmckee} lists $\ell$-torsion frequency data and $c_p$ values for elliptic curves, which is relevant to Conjecture~\ref{galbraithmckee} and the corresponding Lemma~\ref{galbraithmckeeaverage}. Table~\ref{T_ourconjecture} lists similar data for Jacobians of genus~2 curves, see Conjectures~\ref{ourconjecture} and~\ref{imposingweierstrass}, and Lemma~\ref{ourconjectureaverage}.
Table~\ref{T_curveitself1} lists $c_p$ values for the number of points on the curves themselves, related to Conjecture~\ref{curveitselfconj}, while Table~\ref{T_curveitself2} gives
experimental trace distributions of genus $2$ curves modulo $\ell$, see Table~\ref{T_curveitself} above.
Tables~\ref{T_rank1} and~\ref{T_rank2} relate to Theorem~\ref{theoremvladut} and Conjecture~\ref{cyclicitygenus2}, concerning the rank of the Jacobians of curves of genus~1 and~2 (respectively). Finally, Table~\ref{T_ext1} supports Conjecture~\ref{ext1} on the case of extension fields in genus~1.

The data in Tables~\ref{T_galbraithmckee}--\ref{T_rank2} was obtained using the \textsc{smalljac} library \cite{smalljac}, based on the algorithms described in \cite{KedlayaSutherland}. Table~\ref{T_ext1} was obtained using the intrinsic \textsc{Magma} \cite{magma} point counting function.
We conducted our tests by sampling random curves $C$ over finite fields $\FF_p$.
We collected data both using fixed primes $p$, and for all primes in a given interval.
For genus 1 we used $p\approx 10^{12}$ and for genus 2 we used $p\approx 10^6$ (except for Table~\ref{T_curveitself1}) so that in both cases $\#\text{Jac}(C)(\mathbb{F}_p)\approx 10^{12}$.
Each test with a fixed prime used a sample size of approximately $10^6$, while our interval tests used $10^2$ curves for each of at least $10^4$ primes.
In order to maximize the performance of the algorithms used to collect the data, we restricted our tests to curves of the form $y^2=f(x)$, where $f$ is a monic polynomial of degree $2g+1$. Therefore, in genus $2$, our experimental data should be compared to the $\mathcal{H}_5^m$-analogues of the conjectures that
deal with $\mathcal{H}_6$ (which according to Theorem~\ref{maintheoremdeg5} only affects the contribution of $\ell=2$, the necessary
adaptations to which can be made using Table~\ref{T_2torsion}).

%Throughout, we used the offset logarithmic integral $\text{Li}(x)$ to estimate the probability of primality of a uniformly randomly chosen integer from a certain Hasse-Weil interval, rather than actually counting the primes.

\begin{table}[h]
\begin{center}
\small
\begin{tabular}{|l|lccccc|}
\hline
$p$&&$\ell=2$&$\ell=3$&$\ell=5$&$\ell=7$&$c_p$\\\hline
$10^{12}+39$  &observed &0.6654&0.3749&0.2507&0.1664&0.5492\\
          &predicted&0.6667&0.3750&0.2500&0.1667&0.5564\\[3pt]
$10^{12}+61$ &observed &0.6662&0.5003&0.2083&0.1664&0.4686\\
          &predicted&0.6667&0.5000&0.2083&0.1667&0.4646\\[3pt]
$10^{12}+63$ &observed &0.6672&0.3756&0.2503&0.1460&0.5600\\
          &predicted&0.6667&0.3750&0.2500&0.1458&0.5642\\[3pt]
$10^{12}+91$ &observed &0.6660&0.4989&0.2089&0.1454&0.4818\\
          &predicted&0.6667&0.5000&0.2083&0.1458&0.4794\\[3pt]
$[10^{12},10^{12}+4\cdot 10^6]$&observed&0.6666&0.4374&0.2396&0.1631&0.5044\\
     &predicted&0.6667&0.4375&0.2396&0.1632&0.5052\\
\hline
\end{tabular}
\caption{
$\ell$-torsion frequencies and $c_p$ values for $C(\FF_p)$
using random elliptic curves $C:y^2=f(x)$ with $f\in\mathcal{H}_3^\text{m}$.
Sample size is $10^6$ (or $10^2$ for $p$ ranging over the interval $[10^{12},10^{12}+4\cdot 10^6]$).
} \label{T_galbraithmckee}
\end{center}
\end{table}

\begin{table}
\begin{center}
\small
\begin{tabular}{@{}|l|lccccc|@{}}
\hline
$p$&&$\ell=2$&$\ell=3$&$\ell=5$&$\ell=7$&$c_p$\\\hline
$10^6+3$  &observed &0.7991&0.3616&0.2395&0.1628&0.3426\\
          &predicted&0.8000&0.3609&0.2396&0.1632&0.3444\\[3pt]
%$10^6+33$ &observed &0.8007&0.3602&0.2397&0.1633&0.3438\\
%          &predicted&0.8000&0.3609&0.2396&0.1632&0.3472\\[3pt]
$10^6+37$ &observed &0.8000&0.4376&0.2393&0.1626&0.3056\\
          &predicted&0.8000&0.4375&0.2396&0.1632&0.3037\\[3pt]
%$10^6+39$ &observed &0.7993&0.3606&0.2395&0.1636&0.3484\\
%          &predicted&0.8000&0.3609&0.2396&0.1632&0.3464\\[3pt]
$10^6+81$ &observed &0.8001&0.3619&0.2066&0.1632&0.3571\\
          &predicted&0.8000&0.3609&0.2067&0.1632&0.3593\\[3pt]
%$10^6+99$ &observed &0.8004&0.3604&0.2400&0.1632&0.3487\\
%          &predicted&0.8000&0.3609&0.2396&0.1632&0.3472\\[3pt]
%$10^6+117$&observed &0.7999&0.3613&0.2400&0.1629&0.3494\\
%          &predicted&0.8000&0.3609&0.2396&0.1632&0.3464\\[3pt]
$10^6+121$&observed &0.8003&0.4376&0.2059&0.1637&0.3197\\
          &predicted&0.8000&0.4375&0.2067&0.1632&0.3189\\[3pt]
$[10^6,2\cdot 10^6]$&observed&0.8000&0.3992&0.2314&0.1604&0.3285\\
     &predicted&0.8000&0.3992&0.2314&0.1602&0.3290\\
\hline
\end{tabular}
\caption{
$\ell$-torsion frequencies and $c_p$ values for $\text{Jac}(C)(\FF_p)$
using random genus~2 curves $C:y^2=f(x)$ with $f\in\mathcal{H}_5^\text{m}$.
Sample size is $10^6$ (or $10^2$ for $p$ ranging over the interval $[10^6,2\cdot 10^6]$).
} \label{T_ourconjecture}
\end{center}
\end{table}

\begin{table}
\begin{center}
\small
\begin{tabular}{@{}|lcccc|@{}}
\hline
&$10^9+7$   & $10^9+9$  &  $10^9+21$  &  $10^9+33$ \\ \hline
observed      & 1.0162   & 1.0738  &  1.0892   &  1.0945\\
predicted     & 1.0194   & 1.0790  &  1.0865   &  1.0898\\
\hline
\end{tabular}
\caption{
$c_p$ values for the number of points on random genus 2 curves $y^2=f(x)$ with
$f\in\mathcal{H}_5^\text{m}$.  Sample size is $10^6$.
The deviations are larger here due to the shorter intervals (of width approximately $8 \cdot 10^{9/2}$
versus $8\cdot 10^6$ and $4\cdot 10^6$ in Tables \ref{T_galbraithmckee} and \ref{T_ourconjecture} above).
} \label{T_curveitself1}
\end{center}
\end{table}

\begin{table}
\begin{center}
\small
\begin{tabular}{@{}|l|r|lccccc|@{}}
\hline
$p$&$\ell$&&$t\equiv 0$& $t\equiv 1$&$t\equiv 2$&$t\equiv 3$&$t\equiv 4$\\\hline
$10^6+3$  & 2 & observed  & 0.4658 & 0.5342 & & &\\
          &   & predicted & 0.4667 & 0.5333 & & &\\[3pt]
          & 3 & observed  & 0.3598 & 0.3205 & 0.3197 & &\\
          &   & predicted & 0.3594 & 0.3203 & 0.3203 & &\\[3pt]
          & 5 & observed  & 0.2072 & 0.1988 & 0.1978 & 0.1981 & 0.1981\\
          &   & predicted & 0.2067 & 0.1982 & 0.1985 & 0.1985 & 0.1982\\[6pt]
$10^6+37$ & 2 & observed  & 0.4653 & 0.5346 & & &\\
          &   & predicted & 0.4667 & 0.5333 & & &\\[3pt]
          & 3 & observed  & 0.3628 & 0.3185 & 0.3186 & &\\
          &   & predicted & 0.3625 & 0.3188 & 0.3188 & &\\[3pt]
          & 5 & observed  & 0.2070 & 0.1982 & 0.1981 & 0.1983 & 0.1984\\
          &   & predicted & 0.2067 & 0.1985 & 0.1982 & 0.1982 & 0.1985\\[6pt]
$10^6+39$ & 2 & observed  & 0.4667 & 0.5332 & & &\\
          &   & predicted & 0.4667 & 0.5333 & & &\\[3pt]
          & 3 & observed  & 0.3593 & 0.3206 & 0.3202 & &\\
          &   & predicted & 0.3594 & 0.3203 & 0.3203 & &\\[3pt]
          & 5 & observed  & 0.2068 & 0.1978 & 0.1983 & 0.1989 & 0.1982\\
          &   & predicted & 0.2066 & 0.1985 & 0.1983 & 0.1983 & 0.1985\\[6pt]
$[10^6,2\cdot 10^6]$ & 2 & observed  & 0.4669 & 0.5331&  & &\\
          &   & predicted & 0.4667 & 0.5333 & & &\\[3pt]
          & 3 & observed  & 0.3609 & 0.3194 & 0.3197 & &\\
          &   & predicted & 0.3625 & 0.3203 & 0.3203 & &\\[3pt]
          & 5 & observed  & 0.2068 & 0.1982 & 0.1984 & 0.1981 & 0.1985\\
          &   & predicted & 0.2067 & 0.1984 & 0.1984 & 0.1985 & 0.1984\\
\hline
\end{tabular}
\caption{
Trace distributions modulo $\ell$ for random genus 2 curves $y^2=f(x)$ with
$f\in\mathcal{H}_5^\text{m}$.  Sample size is $10^6$ (or $10^2$ for $p$ ranging over the interval $[10^6,2\cdot 10^6]$).
} \label{T_curveitself2}
\end{center}
\end{table}

\begin{table}
\begin{center}
\small
\begin{tabular}{@{}|l|r|llll|@{}}
\hline
$p$&$\ell$&&rank 0& rank 1&rank 2\\\hline
$10^{12}+39$  & 2 & observed  & 0.3346 & 0.4993 & 0.1661\\
          &   & predicted & 0.3333 & 0.5000 & 0.1667 \\[3pt]
          & 3 & observed  & 0.6251 & 0.3334 & 0.0415\\
          &   & predicted & 0.6250 & 0.3333 & 0.0417\\[3pt]
		  & 5 & observed  & 0.7492 & 0.2507 & \\
          &   & predicted & 0.7500 & 0.2500 & \\[3pt]
          & $\infty$ & observed  &         & 0.7988 & 0.2013\\
          &   & predicted & & 0.7980 & 0.2020\\[6pt]
$10^{12}+61$ & 2 & observed  & 0.3338 & 0.4996 & 0.1666\\
          &   & predicted & 0.3333 & 0.5000 & 0.1667\\[3pt]
          & 3 & observed  & 0.4997 & 0.5003 & \\
          &   & predicted & 0.5000 & 0.5000 & \\[3pt]
          & 5 & observed  & 0.7917 & 0.1999 & 0084\\
          &   & predicted & 0.7917 & 0.2000 & 0083\\[3pt]
          & $\infty$ & observed  &         & 0.8263 & 0.1737\\
          &   & predicted & & 0.8264 & 0.1736\\[6pt]
$10^{12}+63$& 2 & observed  & 0.3328 & 0.4995 & 0.1677\\
          &   & predicted & 0.3333 & 0.5000 & 0.1667\\[3pt]
          & 3 & observed  & 0.6244 & 0.3339 & 0.0416\\
          &   & predicted & 0.6250 & 0.3333 & 0.0417\\[3pt]
          & 5 & observed  & 0.7497 & 0.2503 & \\
          &   & predicted & 0.7500 & 0.2500 & \\[3pt]
          & $\infty$ & observed  &  & 0.7953 & 0.2047\\
          &   & predicted &  & 0.7962 & 0.2038\\[6pt]
$[10^{12},2\cdot 10^{12}+4\cdot 10^6]$& 2 & observed  & 0.3334 & 0.4999 & 0.1666\\
          &   & predicted & 0.3333 & 0.5000 & 0.1667\\[3pt]
          & 3 & observed  & 0.5626 & 0.4166 & 0.0208\\
          &   & predicted & 0.5635 & 0.4167 & 0.0208\\[3pt]
          & 5 & observed  & 0.7604 & 0.2375 & 0.0021\\
          &   & predicted & 0.7604 & 0.2375 & 0.0021\\[3pt]
          & $\infty$ & observed  &         & 0.8138 & 0.1862\\
          &   & predicted & & 0.8138 & 0.1862\\
\hline
\end{tabular}
\caption{
Rank frequencies for $C(\FF_p)$ for random elliptic curves $C:y^2=f(x)$ with
$f\in\mathcal{H}_3^\text{m}$.  Sample size is $10^6$ (or $10^2$ for $p$ ranging over the interval $[10^{12},2\cdot 10^{12}]$).
Rows with $\ell=\infty$ indicate maximum $\ell$-rank over all primes $\ell$.
} \label{T_rank1}
\end{center}
\end{table}

\begin{table}
\begin{center}
\small
\begin{tabular}{@{}|l|r|llllll|@{}}
\hline
$p$&$\ell$&&rank 0& rank 1&rank 2&rank 3&rank 4\\\hline
$10^6+3$  & 2 & observed  & 0.200113 & 0.416313 & 0.291528 & 0.083775 & 0.008271\\
          &   & predicted & 0.200000 & 0.416667 & 0.291667 & 0.083333 & 0.008333\\[3pt]
          & 3 & observed  & 0.637964 & 0.320212 & 0.040254 & 0.001548 & 0.000022\\
          &   & predicted & 0.639063 & 0.319444 & 0.039931 & 0.001543 & 0.000019\\[3pt]
		  & 5 & observed  & 0.761095 & 0.236804 & 0.002101 &  & \\
          &   & predicted & 0.760417 & 0.237500 & 0.002083 &  & \\[3pt]
%          & 7 & observed  & 0.836718 & 0.162824 & 0.000458 & 0 & 0\\
%          &   & predicted & 0.836806 & 0.162698 & 0.000496 & 0 & 0\\
          & $\infty$ & observed  &         & 0.589030 & 0.317489 & 0.085188 & 0.008293\\
          &   & predicted &  & 0.589471 & 0.317443 & 0.084733 & 0.008352\\[6pt]
$10^6+81$ & 2 & observed  & 0.200794 & 0.416446 & 0.290857 & 0.083593 & 0.008310\\
          &   & predicted & 0.200000 & 0.416667 & 0.291667 & 0.083333 & 0.008333\\[3pt]
          & 3 & observed  & 0.637636 & 0.320698 & 0.040107 & 0.001533 & 0.000026\\
          &   & predicted & 0.639063 & 0.319444 & 0.039931 & 0.001543 & 0.000019\\[3pt]
          & 5 & observed  & 0.793657 & 0.198090 & 0.008186 & 0.000067 & 0.000000\\
          &   & predicted & 0.793336 & 0.198333 & 0.008264 & 0.000067 & 0.000000\\[3pt]
%          & 7 & observed  & 0.837363 & 0.162134 & 0.000503 & 0 & 0\\
%          &   & predicted & 0.836806 & 0.162698 & 0.000496 & 0 & 0\\
          & $\infty$ & observed  &         & 0.586416 & 0.320192 & 0.085056 & 0.008336\\
          &   & predicted & & 0.585781 & 0.321073 & 0.084794 & 0.008353\\[6pt]
$10^6+133$& 2 & observed  & 0.199300 & 0.416997 & 0.292156 & 0.083233 & 0.008314\\
          &   & predicted & 0.200000 & 0.416667 & 0.291667 & 0.083333 & 0.008333\\[3pt]
          & 3 & observed  & 0.562514 & 0.416732 & 0.020754 & & \\
          &   & predicted & 0.562500 & 0.416667 & 0.020833 & & \\[3pt]
          & 5 & observed  & 0.760019 & 0.237919 & 0.002062 & & \\
          &   & predicted & 0.760417 & 0.237500 & 0.002083 & & \\[3pt]
%          & 7 & observed  & 0.854942 & 0.142071 & 0.002978 & 0.000009 & 0.000000\\
%          &   & predicted & 0.854592 & 0.142432 & 0.002967 & 0.000009 & 0.000000\\
          & $\infty$ & observed  & & 0.600296 & 0.308148 & 0.083242 & 0.008314\\
          &   & predicted & & 0.600635 & 0.307690 & 0.083341 & 0.008333\\[6pt]
$[10^6,2\cdot 10^6]$& 2 & observed  & 0.200039 & 0.416528 & 0.291761 & 0.083320 & 0.008353\\
          &   & predicted & 0.200000 & 0.416667 & 0.291667 & 0.083333 & 0.008333\\[3pt]
          & 3 & observed  & 0.600830 & 0.368047 & 0.030337 & 0.000777 & 0.000009\\
          &   & predicted & 0.600781 & 0.368056 & 0.030382 & 0.000772 & 0.000010\\[3pt]
          & 5 & observed  & 0.768609 & 0.227739 & 0.003637 & 0.000016 & 0.000000\\
          &   & predicted & 0.768647 & 0.227708 & 0.003629 & 0.000017 & 0.000000\\[3pt]
%          & 7 & observed  & 0.839643 & 0.159431 & 0.000925 & 0.000001 & 0.000000\\
%          &   & predicted & 0.839770 & 0.159321 & 0.000908 & 0.000001 & 0.000000\\
          & $\infty$ & observed  &         & 0.594471 & 0.313125 & 0.084043 & 0.008362\\
          &   & predicted & & 0.594567 & 0.313040 & 0.084050 & 0.008343\\
\hline
\end{tabular}
\caption{
Rank frequencies for $\text{Jac}(C)(\FF_p)$ for random genus 2 curves $C:y^2=f(x)$ with
$f\in\mathcal{H}_5^\text{m}$.  Sample size is $10^6$ (or $10^2$ $p$ ranging over the interval $[10^6,2\cdot 10^6]$).
Rows with $\ell=\infty$ indicate maximum $\ell$-rank over all primes $\ell$.
} \label{T_rank2}
\end{center}
\end{table}

\begin{table}
\begin{center}
\small
\begin{tabular}{|l|lcccccc|} \hline
$p^k$   &                 & $\ell=2$   & $\ell=3$   & $\ell=5$   & $\ell=7$   & $\ell=11$   & $c_k$    \\\hline

$3^{26}$  & \text{observed} & $0.6669$   & $0.4999$   & $0.2501$   & $0.1666$   & $0.1000$    & $0.4387$ \\
   & \text{predicted}& $0.6667$   & $0.5000$   & $0.2500$   & $0.1667$   & $0.1000$    & $0.4401$ \\[3pt]

$5^{18}$  & \text{observed} & $0.6667$   & $0.3748$   & $0.2501$   & $0.1458$   & $0.1000$    & $0.5659$ \\
 & \text{predicted}& $0.6667$   & $0.3750$   & $0.2500$   & $0.1458$   & $0.1000$    & $0.5662$ \\[3pt]

$7^{15}$  & \text{observed} & $0.6667$   & $0.3751$   & $0.2499$   & $0.1667$   & $0.1001$    & $0.5541$ \\
 & \text{predicted}& $0.6667$   & $0.3750$   & $0.2500$   & $0.1667$   & $0.1000$    & $0.5523$ \\[3pt]

$11^{12}$  & \text{observed} & $0.6665$   & $0.3749$   & $0.2083$   & $0.1457$   & $0.1002$    & $0.6020$ \\
 & \text{predicted}& $0.6667$   & $0.3750$   & $0.2083$   & $0.1458$   & $0.1000$    & $0.6015$ \\\hline
\end{tabular}
\caption{$\ell$-torsion frequencies and $c_k$ values
for $C(\mathbb{F}_{p^k})$ using random elliptic curves $C:y^2=f(x)$
with $f\in\mathcal{H}_3^\text{m}$. Sample
size is $10^7$.} \label{T_ext1}
\end{center}
\end{table}

\clearpage

\section{Acknowledgements}

We would like to thank Steven Galbraith for proposing this research, and Jeff Achter, Jason Fulman, Frans Oort, Bjorn Poonen,
Alessandra Rigato, Igor Shparlinski, Marco Streng, and the anonymous referee for some helpful comments. The first author is grateful to the Massachusetts Institute of Technology for its hospitality, and to F.W.O.-Vlaanderen for its financial support.

\noindent \emph{Katholieke Universiteit Leuven}\\
\noindent \emph{Departement Wiskunde}\\
\noindent \emph{Celestijnenlaan 200B, 3001 Leuven (Heverlee), Belgium}\\
\noindent \verb"wouter.castryck@gmail.com"\\

\noindent \emph{Yale University}\\
\noindent \emph{Mathematics Department}\\
\noindent \emph{P.O.\ Box 208283, New Haven, CT 06520-8283, USA}\\
\noindent \verb"amanda.folsom@yale.edu"\\

\noindent \emph{Katholieke Universiteit Leuven}\\
\noindent \emph{Departement Wiskunde}\\
\noindent \emph{Celestijnenlaan 200B, 3001 Leuven (Heverlee), Belgium}\\
and\\
\noindent \emph{Universit\'{e} Libre de Bruxelles}\\
\noindent \emph{D\'{e}partement de Math\'{e}matique}\\
\noindent \emph{Boulevard du Triomphe, 1050 Brussels, Belgium}\\
\noindent \verb"hendrik.hubrechts@wis.kuleuven.be"\\

\noindent \emph{Massachusetts Institute of Technology}\\
\noindent \emph{Department of Mathematics}\\
\noindent \emph{77 Massachusetts Avenue, Cambridge, MA 02139-4307, USA}\\
\noindent \verb"drew@math.mit.edu"\\
\end{document}